\documentclass[11pt,reqno]{amsart}

\usepackage{xfrac}
\usepackage{amsmath,amsthm,amsfonts,epsfig,color,esint,mathtools,amssymb}
\usepackage[latin1]{inputenc}
\usepackage{hyperref}

\newtheorem{theorem}{Theorem}[section]
\newtheorem{lemma}[theorem]{Lemma}
\newtheorem{proposition}[theorem]{Proposition}
\newtheorem{definition}[theorem]{Definition}

\newtheorem{conjecture}[theorem]{Conjecture}

\newtheorem{problem}[theorem]{Problem}

\renewcommand*{\div}{\ensuremath{\mathrm{div\,}}}

\newcommand\tr{{\rm tr}\,}
\def\deg{\textrm{deg}\,}

\def\eps{\varepsilon}

\def\S{\mathcal{S}}

\def\R{\mathbb{R}}
\def\T{\mathbb{T}}
\def\Z{\mathbb{Z}}
\def\N{\mathbb{N}}
\def\C{\mathbb{C}} 
\def\Id{{\rm Id}}

\begin{document}



\title{High-dimensionality and h-principle in PDE}

\author{Camillo De Lellis}
\address{Institut f\"ur Mathematik, Universit\"at Z\"urich, 
CH-8057 Z\"urich}
\email{camillo.delellis@math.uzh.ch}
\thanks{The research of C.D.L. has been supported by the grant $200021\_159403$ of the Swiss National Foundation.}

\author{L\'aszl\'o Sz\'ekelyhidi Jr.}
\address{Mathematisches Institut, Universit\"at Leipzig, 
D-04009 Leipzig}
\email{laszlo.szekelyhidi@math.uni-leipzig.de}
\thanks{L.Sz. gratefully acknowledges the support of ERC Grant Agreement No.~277993 }

\maketitle

\tableofcontents

\section{Introduction}

In this same Volume of the Bulletin of the AMS M.~Gromov wrote a highly inspiring and visionary article on the legacy of John Nash's papers on isometric embeddings. There is surely no better choice of author for such a task - indeed, it was Gromov who realized that Nash's works on isometric embeddings not only solved the existence problem as it was formulated at the time, but actually opened the door to a completely new type of mathematics reaching far beyond differential geometry. We cannot resist to quote him directly: {\it What Nash discovered is not any part of Riemannian geometry, neither has it much (if anything at all) to do with classical PDE.}  

In his article Gromov paints a picture of the ``New Land'' discovered by Nash with exceptional clarity, breath and depth. One of the key aspects of Nash's theorems emphasized in the article is the {\it high dimensionality} (``Infinite dimensional representations''). Whilst in the case of the smooth embedding theorem of Nash \cite{Nash1956} the high-dimensionality is of geometric nature, in the $C^1$ embedding theorem of Nash-Kuiper \cite{Nash1954,Kuiper} it is rather of analytic nature. In this note we would like to present ``an analysts' point of view'' on the latter theorem and in particular highlight the very close connection to turbulence -- a paradigm example of a high-dimensional phenomenon! Our aim is to explain recent applications of Nash's ideas in connection with the incompressible Euler equations and Onsager's famous conjecture on the energy dissipation in 3D turbulence.

 \section{Isometric Embeddings, Nash, and Gromov's h-principle}

In his book \cite{Gromov86} Gromov developed convex integration, a far reaching generalization of the perturbation technique of Nash 1954/1956. The most common application of convex integration is to provide solutions to a certain generic class of partial differential relations, consisting of a global topological condition and a {\it differential inequality} representing the non-singularity of some geometric quantity, i.e. the non-vanishing of some function of the derivatives. Examples include the Smale's sphere inversion and the existence of $n$ linearly independent non-vanishing divergence-free vector fields (see the book of Gromov \cite{Gromov86} and his current article). In such problems if a solution exists at all, then -- obviously -- there exist infinitely many solutions: indeed, the solution space is \emph{open} in an appropriate function space. It is a curious fact of life that finding a solution becomes much more difficult if there is no uniqueness (even locally!), because then, while looking for a solution, there is no way to characterize it, or at least to formulate a clear preference. In simplified terms convex integration produces a large family of local perturbations which keep the topological condition whilst achieving the required non-vanishing. Of course the situation is in reality more complicated: for instance, one can easily ensure by a local perturbation that the derivative of a function on the unit circle $f:{\mathbb S}^1\to\R$ is not zero in any given small neighborhood; however it will always have zero derivative \emph{somewhere}. Indeed, the main issue is to understand how the global topology affects the local differential structure. For situations where the topology ``wins'' over the local geometry, Gromov introduced the term ``h-principle''.  

\subsection{A classical problem in Differential Geometry}
The problem of isometric embeddings of Riemannian manifolds consists of a global topological condition (being an embedding) and a system of partial differential equations (being an isometry), so it seems at first glance to be  a completely different kind of problem.
For concreteness let us consider  a smooth $n$-dimensional Riemannian manifold $(\Sigma^n, g)$. A continuous map $u : \Sigma \to \mathbb R^N$ is {\em isometric} if it preserves the length of curves, namely if
\begin{equation}\label{e:true_isometry}
\ell_g (\gamma) = \ell_e (u\circ \gamma) \qquad \mbox{for any $C^1$ curve $\gamma\subset \Sigma$,}
\end{equation}
where $\ell_g (\gamma)$ denotes the length of $\gamma$ with respect to the metric $g$:
\begin{equation}\label{e:Riem_length}
\ell_g (\gamma) = \int \sqrt{ g (\gamma (t)) [\dot\gamma (t), \dot\gamma (t)]}\, dt\, .
\end{equation}
 As customary, in local coordinates we can express the metric tensor $g$ as\footnote{Here and in the rest of this note we follow Einstein's summation convention.}
\[
g = g_{ij} dx_i \otimes dx_j\, ,
\]
For a general $C^1$ map $v$ the tensor $(\partial_i v \cdot \partial_j v)\,  dx_i \otimes dx_j$ is usually called the ``pull-back'' of the euclidean metric and thus it is customary to denote it by $v^\sharp e$; in general $v^\sharp e$ is only positive semidefinite and it is positive definite (hence a metric) if and only if $v$ is an immersion. With this notation at hand we can rewrite the condition \eqref{e:true_isometry} for $C^1$ maps as $u^\sharp e = g$ and such identity alone guarantees that $u$ is an immersion. Then, if $u$ is $C^1$, \eqref{e:true_isometry} is equivalent to a system of partial differential equations, which in local coordinates takes the following form:
\begin{equation}\label{e:equations}
\partial_iu\cdot\partial_ju=g_{ij}\, .
\end{equation} 

The existence of isometric immersions (resp.~embeddings) of Riemannian manifolds into some Euclidean space is a classical problem, explicitly formulated for the first time by Schl\"afli, see \cite{Schlaefli}. Clearly, if the dimension of $\Sigma$ is $n$, \eqref{e:equations} consists of $s_n:=\frac{n(n+1)}{2}$ equations in $N$ unknowns. 
A reasonable guess would therefore be that the system is solvable, at least locally,
when $N= s_n$: this was in fact what Schl\"afli conjectured in his note. 

In the first half of the twentieth century Janet
\cite{Janet}, Cartan \cite{Cartan} and Burstin \cite{Burstin} had proved the existence of local isometric embeddings in the case of analytic metrics, precisely when $N= s_n$. For the very particular case of $2$-dimensional spheres endowed with metrics of positive Gauss curvature, Weyl in \cite{Weyl} had raised the question of the existence of (global!) isometric embeddings in $\mathbb R^3$. The Weyl's problem was solved by Lewy in \cite{Lewy} for analytic metrics and Louis Nirenberg settled the case of smooth metrics in his PhD thesis in 1949 (in fact Nirenberg's Theorem requires $C^4$ regularity of the metric tensor, see \cite{NirenbergPhD} and \cite{Nirenberg}); a different proof was given independently by Pogorelov \cite{Pogorelov1951} around the same time, building upon the work of Alexandrov \cite{Alexandrov1948} (see also \cite{Pogorelov73}). 
Moreover in the case of the Weyl's problem it was proved by Herglotz and Cohn-Vossen, already before the work of Lewy, that $C^2$ immersions are uniquely determined up to rigid motions, cf. \cite{CohnVossen,Herglotz} and see also \cite{Spivak5} for a thorough discussion. 
Incidentally, the linearized (infinitesimal) rigidity in the Weyl problem, due to Blaschke \cite{Blaschke}, was of crucial importance in Nirenberg's existence proof -- a nice example of how uniqueness leads to existence.

\subsection{The paradox of Nash} Nash started working at Schl\"afli's general question, which was considered a formidable problem, shortly after his PhD, apparently because of a bet with a colleague at the MIT department, where he had just moved as a young faculty, cf. \cite{Nasar}. As John Milnor wrote recently ``Nash was never a reasonable person'' and indeed, although at the time everything indicated that the solvability of \eqref{e:equations} needs a high dimensional target, in his 1954 note \cite{Nash1954} Nash astonished the geometry world and proved that the only true obstructions to the existence of isometric immersions are topological. As soon as $N\geq n+1$ and there are no such obstructions, then there are in fact plenty of such immersions.
 
Nash gave a proof of this statement for $N\geq n+2$ and just remarked that a similar one could be proved for $N\geq n+1$; the details were then given in two subsequent notes by Kuiper, \cite{Kuiper}. For this reason the resulting theorem is called nowadays the Nash-Kuiper Theorem on $C^1$ isometric embeddings. In order to state it, we follow Nash and introduce first the notion of ``short maps'', namely maps which decrease lengths.

\begin{definition}\label{d:short}
Let $(\Sigma, g)$ be a Riemannian manifold. 
An immersion $v: \Sigma \to \mathbb R^N$ is short if we have the inequality $h:= v^\sharp e \leq g$ in the sense of quadratic forms,
i.e. in local coordinates we have $h_{ij} w^i w^j \leq g_{ij} w^i w^j$ for any tangent vector $w$.
If the strict inequality $<$ holds we then say that $v$ is {\em strictly} short.
\end{definition}
The Nash-Kuiper Theorem is then the following

\begin{theorem}\label{t:main_C1_1}
Let $(\Sigma, g)$ be a smooth closed $n$-dimensional Riemannian manifold and $v: \Sigma \to \mathbb R^N$ a $C^\infty$ short immersion with $N\geq n+1$. Then, for any $\varepsilon >0$ there exists a $C^1$ isometric immersion $u: \Sigma\to \mathbb R^N$ such that $\|u-v\|_{C^0} \leq \varepsilon$. If $v$ is, in addition, an embedding, then $u$ can be assumed to be an embedding as well. 
\end{theorem}

This theorem shows -- and it was Gromov who understood the deep implications of this interpretation -- that the system \eqref{e:equations} of non-linear partial differential equations is sufficiently ``soft'' so that in a certain sense it behaves more like a differential inequality.  In particular, although the set of isometries obviously cannot form an open set in the space of $C^1$ maps, it is $C^0$-dense in the open set of strictly short maps. This type of abundance of solutions is a central aspect of Gromov's $h$-principle.
In addition note that when $\Sigma$ is a smooth closed manifold we can make any immersion $v: \Sigma \to \mathbb R^N$ short by simply multiplying it by a small positive constant. Hence Theorem \ref{t:main_C1_1} reduces the existence of isometries (resp. isometric embeddings) to that of immersions (resp. embeddings), which is guaranteed by the classical Theorem of Whitney in a codimension which is rather low compared to the codimension in Schl\"afli's conjecture. 

\section{Soft PDEs and thresholds}
\subsection{Relaxation}\label{s:relaxation1}
A good characterization of nonlinear differential structures which are soft is still missing, although partial answers based on L.~Tartar's formalism, compensated compactness and relaxation exist, see for instance \cite{KirMulSver,DS-Survey}. In order to explain the basic idea of this approach let us again look at the system of partial differential equations \eqref{e:equations} with some fixed smooth $g$, and consider a sequence of (smooth) solutions $\{u^k\}_k$, $u^k: {\mathbb S}^2 \to \R^3$. Then the sequence of derivatives $|\partial_i u^k|^2$ is uniformly bounded, hence by the Arzel\`a-Ascoli theorem there exists a subsequence $u^\kappa$ converging uniformly to some limit map $u$. The limit $u$ must be Lipschitz and an interesting question is whether $u$ is still a solution, i.e. isometric. This would follow from some better convergence, for instance in the $C^1$ category.
If the metric $g$ has positive curvature and the maps $u^k$ are sufficiently smooth, their images will be convex surfaces: this, loosely speaking, amounts to some useful information about second derivatives which will improve the convergence of the subsequence $u^\kappa$ and result in a limit $u$ with convex image.

If instead we only assume that the sequence $u^k$ consists of approximate solutions, for instance in the sense that
\[
\partial_iu^k\cdot \partial_ju^k-g_{ij}\to 0\textrm{ uniformly,}
\]
then even if $g$ has positive curvature and the $u^k$ are smooth, their images will not necessarily be convex, as we can (nowadays!) infer from the Nash-Kuiper Theorem.
 
Let us nonetheless see what we can conclude about the limit $u$. Consider a smooth curve $\gamma\subset {\mathbb S}^2$. Then $u^k \circ \gamma$ is a $C^1$ Euclidean curve and our assumption implies  
\begin{equation}\label{e:lengths}
\ell(u^k \circ \gamma) \to \ell (\gamma).
\end{equation}
On the other hand the curves $u^k \circ \gamma$ converge uniformly to the Lipschitz curve $u\circ \gamma$ and it is well-known that under such type of convergence the length might shrink but cannot increase. We conclude that
\begin{equation}\label{e:l_short}
\ell (u \circ \gamma) \leq \ell(\gamma)\, ,
\end{equation}
in other words the map $u$ is {\it short}.
Recall that, by Rademacher's theorem, $u$ is differentiable almost everywhere: it is a simple exercise to see that,
when \eqref{e:l_short} holds for every Lipschitz curve $\gamma$, then
\begin{equation}\label{e:short}
\partial_i u \cdot \partial_j u \leq g_{ij} \qquad \mbox{a.e.,}
\end{equation}
in the sense of quadratic forms. Thus, loosely speaking, one possible interpretation of  Theorem \ref{t:main_C1_1} is that the system of partial differential inequalities \eqref{e:short} is the ``relaxation'' of \eqref{e:equations} with respect to the $C^0$ topology.

\subsection{Thresholds} 
Now, for $C^2$ isometric immersions there are higher order constraints, most notably the Theorema Egregium of Gauss. This is in fact one crucial ingredient in the proof of rigidity for the Weyl problem. In particular such rigidity 
implies that any $C^2$ isometric immersion of the standard sphere in $\mathbb R^3$ must map it to the boundary of some ball of radius $1$. On the contrary the Nash-Kuiper theorem implies the existence of $C^1$ isometric embeddings which crumple the standard sphere into an arbitrarily small region of the $3$-dimensional space. 

This implies a counterintuitive dichotomy between ``rough'' and ``smooth''  (i.e. below and above $C^2$) solutions of \eqref{e:equations} in low codimension.
An interesting open question, which will be explored further in this note, is whether there exists a ``threshold regularity'' which distinguishes between these two phenomena. A particular case of this question is the following

\begin{problem}\label{p:2d-Hoelder}
Let $N = 3 = n+1$.
Is there a threshold $\theta_0\in ]0,1[$ such that:
\begin{itemize}
\item $C^{1, \theta}$ solutions of the Weyl problem are rigid for $\theta> \theta_0$;
\item the Nash-Kuiper Theorem holds for $C^{1,\theta}$ immersions when $\theta < \theta_0$?
\end{itemize}
\end{problem}

Indeed the question of which regularity one might reach with Nash's 1954 scheme turns out to be an interesting and difficult open problem, with ramifications beyond the isometric embedding problem. 
In particular, as we have pointed out recently, a celebrated conjecture of Lars Onsager in the theory of fully developed turbulence shares many similarities with Problem \ref{p:2d-Hoelder} and can be approached with an iteration which is similar to Nash's 1954 scheme, see \cite{DS-Annals, DS-Inv, BDIS, Sz-ICM}. Although we will discuss its context and the precise definitions later, we state here the Onsager's conjecture so that the reader could appreciate the formal analogy with Problem \ref{p:2d-Hoelder}.

\begin{conjecture}\label{c:Onsager}
Consider periodic $3$-dimensional weak solutions of the incompressible Euler equations, where the velocity $v$ satisfies the uniform H\"older condition
\begin{equation}\label{e:hoelderest}
|v(x,t)-v(x',t)|\leq C|x-x'|^\theta,
\end{equation}
for constants $C$ and $\theta$ independent of $x,x'$ and $t$.
\begin{enumerate}
\item[(a)] If $\theta>\frac{1}{3}$, then the total kinetic energy of $v$ is constant;
\item[(b)] For any $\theta<\frac{1}{3}$ there are $v$ for which it is {\em not} constant. 
\end{enumerate}
\end{conjecture}

Of course Problem \ref{p:2d-Hoelder} deals with a stronger property, namely the rigidity (and thus ``uniqueness'') of the solution. A more stringent analog in the case of the Euler equations would then claim an appropriate uniqueness result for the velocity $v$, for instance for the corresponding Cauchy problem). On the other hand, as already mentioned (and will be explained briefly in Section \ref{s:speculations}), a crucial point in Problem \ref{p:2d-Hoelder} is whether a suitable version of Gauss' Theorema Egregium holds or not at low regularity. We can regard Gauss' Theorema Egregium as an additional identity valid for sufficiently regular solutions of \eqref{e:equations}, pretty much as the conservation law for the energy is an additional identity that sufficiently regular solutions of the Euler equations must fulfill. 

\medskip

In the rest of this note: 
\begin{itemize}
\item we will review Nash's approach to Theorem \ref{t:main_C1_1}, highlighting its ``nonlinear flavor'' (cf. Section \ref{s:Nash-scheme}); 
\item we will give a survey on the state of the art for Problem \ref{p:2d-Hoelder} and related questions (cf. Section \ref{s:Borisov}); 
\item we will give a survey on the most recent results on the Onsager's conjecture (cf. Sections \ref{s:Euler}); 
\item we will discuss an analog of Nash's iteration which produces counterintuitive continuous solutions of the Euler
equations (cf. Section \ref{s:Onsager}); 
\item we will explain how suitable adjustments in the latter iteration leads to H\"older solutions, cf. Section \ref{s:transport}, and to a 
related $h$-principle statement, cf. Section \ref{s:h-principle};
\item we will point out further directions and related open questions (cf. Section \ref{s:speculations}). 
\end{itemize}

\section{Nash's 1954 scheme}\label{s:Nash-scheme}

In his subsequent celebrated 1956 note on the topic (see \cite{Nash1956}) Nash turned his attention to {\em more regular} isometric immersions (resp. embeddings). In particular he proved their existence if the dimension $N$ is sufficiently high, in fact larger than what Schl\"afli conjectured.  If the Nash-Kuiper Theorem could be regarded as a curiosity, the 1956 paper gave a final proof that the abstract worlds of Riemann coincide completely with the usual Euclidean submanifolds. It is well known that 
the impact of Nash's second work goes way beyond its specific application to the isometric embedding problem: his celebrated strategy to treat ``hard implicit function theorems'' has had a profound influence in analysis, in mathematical physics and in geometry in the subsequent 60 years. 
The impact of his first note has been, comparatively, much more modest.

And yet we wish to point out that even the 1954 paper reaches far beyond differential geometry. Its approach to construct solutions of \eqref{e:equations} can be regarded as a fully nonlinear iteration scheme that is highly original and might be applied to other partial differential equations. In classical perturbation methods for nonlinear equations the linearization plays the key role: in this sense the 1956 scheme of Nash is no exception. In contrast, in the 1954 scheme the leading order term is quadratic in the perturbation and the linearization becomes negligible. 
That scheme is thus genuinely infinite dimensional and it is not entirely surprising that it leads to highly irregular solutions.

\subsection{Stage} Let us start reviewing the main ideas in Nash's proof of Theorem \ref{t:main_C1_1}. 
Let $(\Sigma,g)$ be a smooth, closed $n$-dimensional Riemannian manifold.
The proof of Theorem \ref{t:main_C1_1} relies on the following proposition, which Nash calls ``a stage'', 
cf.~\cite[Page 391]{Nash1954}:

\begin{proposition}\label{p:iter_C1}
Let $u: \Sigma \to \mathbb R^N$ be a smooth strictly short immersion. For any $\delta>0$ there
exists a smooth strictly short immersion $\tilde u: \Sigma \to \mathbb R^N$ such that\footnote{As usual, for maps 
$u:\Sigma\to\R^N$ we define the $C^1$ seminorm $[u]_{1} = \|Du\|_0$ using an atlas of smooth charts $\{U_\alpha\}$ on $\Sigma$. Moreover, for any symmetric $(0,2)$ tensor $h$ on $\Sigma$ we denote by $\|h\|_0$ the supremum of the Hilbert Schmidt norm of the matrices $h_{ij} (p)$ for 
$p\in \Sigma$. }
\begin{align}
\|u-\tilde u\|_0 & \leq \delta\, , \label{e:iter_1}\\ 
\|g - \tilde u^\sharp e\|_0 & \leq \delta\, ,\label{e:iter_2}\\
[u-\tilde u]_1 & \leq C \|g - u^\sharp e\|_0^{\sfrac12}\, ,\label{e:iter_3}
\end{align} 
for a constant $C$ which depends only upon $\Sigma$. If $u$ is injective, then $\tilde u$ is also injective.
\end{proposition}

Even without estimate \eqref{e:iter_3} this proposition is quite powerful. It says that the set of \emph{almost isometries} is dense (in the uniform topology) in the set of short immersions - a first hint at the type of relaxation statement and underlying h-principle explained in the previous sections. In fact, this type of global approximation statement is key not only in the proof of the Nash-Kuiper theorem on $C^1$ isometries, but also in the proof of Nash's theorem on $C^\infty$-isometries \cite{Nash1956}.

\subsection{Steps and spirals}\label{s:spirals}
 The main idea behind Proposition \ref{p:iter_C1} is the following simple perturbation step:
Let $u: \Sigma \to \mathbb R^N$ be a smooth immersion and let $U\subset\Sigma$ be a single chart with local coordinates $(x_1,\dots,x_n)$. Assuming, as Nash does\footnote{The extension to $N=n+1$ is contained in the papers of Kuiper \cite{Kuiper}. The main difference is the form of the perturbation; instead of a spiral as in \eqref{e:spiral} one needs to use a corrugation, which cannot be written down quite so explicitly. We refer the interested reader also to \cite{ConDS,Sz-Lecturenotes}.}, that $N\geq n+2$, there exist two linearly independent
unit normal vectors $\zeta, \eta$ to $u(U)$, i.e.  $\zeta,\eta: U\to\R^N$ such that for any $i=1,\dots,n$
\begin{equation}\label{e:normalvectors}
|\zeta|=|\eta|=1,\quad \zeta\cdot\eta=0,\textrm{ and }\,\partial_iu\cdot \zeta=\partial_iu\cdot \eta=0.
\end{equation}
Next, let $\xi$ be a unit vector in $\R^n$ and set
\begin{equation}\label{e:spiral}
\tilde u(x):=u(x)+\frac{a(x)}{\lambda}\biggl(\sin(\lambda x\cdot\xi)\zeta(x)+\cos(\lambda x\cdot \xi)\eta(x)\biggr)
\end{equation}
for some amplitude $a=a(x)$ and frequency $\lambda\gg 1$. One directly calculates:
\begin{equation}\label{e:derivative}
\partial_i\tilde u(x)=\partial_iu(x)+a(x)\biggl(\cos(\lambda x\cdot\xi)\zeta\xi_i-\sin(\lambda x\cdot\xi)\eta\xi_i\biggr)+O\left(\frac{1}{\lambda}\right),
\end{equation}
so that, because of \eqref{e:normalvectors}
\begin{equation}\label{e:newmetric}
\partial_i\tilde u\cdot \partial_j\tilde u=\partial_iu\cdot \partial_ju+a^2(x)\xi_i\xi_j+O\left(\frac{1}{\lambda}\right).
\end{equation}
In other words, the spiral perturbation in \eqref{e:spiral} leads to a new map $\tilde u$, whose induced metric, given by  \eqref{e:newmetric}
is -- up to an error of size $\lambda^{-1}$ -- \emph{increased} in the direction of $\xi$ by an amount $a^2$ and is \emph{essentially not changed} in orthogonal directions. By a suitable (large) choice of $\lambda$ we then achieve:
\begin{lemma}\label{l:step}
Let $u:\Sigma\to\R^N$ be a smooth immersion. Let $U\subset\Sigma$ be an open subset of $\Sigma$ contained in a single chart, $a\in C_c^{\infty}(U)$ a smooth function with compact support and $\xi\in \R^n$ a unit vector. For any $\delta>0$ there exists a smooth immersion $\tilde u: \Sigma \to \R^N$ such that:
\begin{align}
\|\tilde u-u\|_0 & \leq \delta\, , \label{e:step_1}\\ 
\|\partial_i\tilde u\cdot \partial_j\tilde u-\partial_iu\cdot \partial_ju-a^2\xi_i\xi_j\|_0 & \leq \delta \, ,\label{e:step_2}\\
[\tilde u - u]_1 & \leq C \|a\|_0\, ,\label{e:step_3}
\end{align} 
for a dimensional constant $C$. 
\end{lemma}

Now, let us assume in addition that $u:\Sigma\to\R^N$ is strictly short. This amounts to the condition that the ``metric error'' 
$h:=g-u^\sharp e$ is positive definite, i.e. it is also a metric on $\Sigma$. Then the implementation of Lemma \ref{l:step} in the proof of Proposition \ref{p:iter_C1} depends ultimately on being able to decompose an arbitrary metric $h$ on $\Sigma$ in a finite sum\footnote{In the case of non-compact manifolds this will indeed be a locally finite sum.} as
\begin{equation}\label{e:decomposition}
h=\sum_\alpha a_{\alpha}^2(d\psi_\alpha)^2,
\end{equation}
where each $a_{\alpha}$ is compactly supported in a single chart $U_\beta\subset\Sigma$ and in local coordinates $\psi_{\alpha}(x)=\sum_{i=1}^n\xi_i^{\alpha}x_i$ for some unit vectors $\xi^\alpha\in\R^n$. Gromov calls this the Kuratowski-Weyl-Nash decomposition, and we refer the reader to his article for interesting generalizations and open questions.

For each term in this decomposition we can apply Lemma \ref{l:step} and use the obvious estimate $\|a_{\alpha}\|_0\leq \|h\|_0$ to obtain a (finite) sequence of corrections $u_0=u,\,u_1,\dots,u_m$, where $m$ is the number of terms in the sum. The final immersion $\tilde u:=u_m$ then satisfies the conclusions of Proposition \ref{p:iter_C1}. 

\subsection{Iteration and convergence}\label{s:Nashiteration}

It is not difficult to prove Theorem \ref{t:main_C1_1} from Proposition \ref{p:iter_C1}, at least for the case of immersions, by a simple iteration\footnote{For embeddings we need an additional argument and refer the reader to \cite{D-Abelvolume,DS-Kohn,Sz-Lecturenotes} for details.}. 

However, it turns out that a restricted version of Proposition \ref{p:iter_C1} already suffices for iteratively removing the error, once we have a ``sufficiently good'' first approximation. In order to simplify the discussion, let us restrict from now on to a single chart; in other words, we assume that $U\subset\R^n$ is a bounded simply connected domain, $g=(g_{ij})$ is a smooth metric (i.e. positive definite form) on $U$ and $u:(\overline{U},g)\to\R^N$ is a smooth strictly short immersion.  

By assumption $(g-Du^TDu)(x)$ is positive definite on $\overline{U}$.  Therefore there exists $\gamma>0$ so that $g-Du^TDu-2\gamma\Id$ is positive definite. In particular $u$ is also a strictly short immersion of the manifold $(U,\tilde g)$, where $\tilde g:=g-\gamma\Id$. Applying Proposition \ref{p:iter_C1} once to $u:(U,\tilde g)\to\R^N$ with $\delta>0$, we obtain a new smooth immersion $u_0$ such that 
$\|Du_0^TDu_0-\tilde g\|_0\leq \delta$. Then $g-Du_0^TDu_0=\gamma \Id+O(\delta)$, and by choosing $\delta>0$ sufficiently small, we may therefore ensure that the new metric error satisfies
\[
(g-Du_0^TDu_0)(x)\in \mathcal{C}_{\delta/\gamma}\quad\textrm{ for all }x,
\]
where
\begin{equation}\label{e:cone}
\mathcal{C}_r:=\left\{A\in {\rm Sym}_{n\times n}:\,\left|\frac{A}{\tfrac{1}{n}|\tr A|}-\Id\right|<r\right\}.
\end{equation}
Geometrically $\mathcal{C}_r$ is a convex cone of positive-definite matrices with opening ``angle'' $r$ centered around the half-line $\{\lambda \Id:\lambda>0\}$. The advantage of introducing the cone $\mathcal{C}_r$ is that it allows us to localize the decomposition \eqref{e:decomposition} in the space of metrics, resulting in a minimal decomposition. This is based on the following elementary linear algebra lemma:

\begin{lemma}\label{l:decomp1}
There exists a dimensional constant $r_0(n)>0$ and $s_n= \frac{n(n+1)}{2}$ unit vectors $\xi^k \in \mathbb R^n$ with the following property. Any matrix $A\in \mathcal{C}_{r_0}$ can be written in a unique way as a \emph{positive} linear combination 
\begin{equation}\label{e:linear_comb}
A = \sum_{k=1}^{s_n} \mu_k^2 (A) \xi^k\otimes \xi^k,
\end{equation} 
where the $\mu_k$ are smooth positive $1/2$-homogeneous functions on $\mathcal{C}_{r_0}$.
\end{lemma}

In other words the set of rank-one semidefinite matrices $\{\xi^k\otimes \xi^k\}$ generates a convex cone of positive semidefinite matrices, which contains $\mathcal{C}_{r_0}$. Since the number $s_n$ is the dimension of the space of symmetric matrices, it is clearly the minimal number for which the decomposition of Lemma \ref{l:decomp1} can be valid in $\mathcal{C}_{r_0}$. A similar decomposition to \eqref{e:linear_comb}, which is valid for all positive definite $A$, can also be proved using a locally finite partition of unity in the space of positive definite matrices (this is contained in the paper of Nash \cite{Nash1954}, see also \cite{D-Abelvolume,Sz-Lecturenotes}), although then the sum in \eqref{e:linear_comb} is only locally finite and the number of non-vanishing terms is significantly larger than $s_n$. Such a decomposition has also proved useful in other contexts, see \cite[Lemma 17.13]{GT} and \cite{MW}.

Next, set $\delta_q=\eps 2^{-q}$ and define for all $q\in\N$
\begin{equation}\label{e:metric_q}
g_q:=g-\delta_{q}\Id\, .
\end{equation}
We construct inductively a sequence of smooth immersions 
\[
u_q:(U,g_q)\to\R^N
\] 
with metric error 
\begin{equation}\label{e:metricerror_q}
\|g_q-Du_q^TDu_q\|_0\leq c_0\delta_{q+1},
\end{equation}
where the dimensional constant $c_0<1$ will be chosen later. Note that here we do not require $u_q$ to be short with respect to the metric $g_q$, but obviously it will be strictly short with respect to the metric $g$. Set 
\begin{equation}\label{e:newmetricerror}
h_q:=g_{q+1}-Du_q^TDu_q.
\end{equation}
It is easy to check, that $h_q(x)\in \mathcal{C}_{r_0}$ for all $x$, provided $c_0$ is sufficiently small (depending only on $r_0$). Therefore we can define the amplitudes
\begin{equation}\label{e:a_k_explicit}
a_{q,k}(x) := \mu_k \bigl(h_q(x)\bigr)
\end{equation}
and obtain from Lemma \ref{l:decomp1} 
\begin{equation}\label{e:decomp_local}
h_q (x) = \sum_{k=1}^{s_n} a_{q,k}^2 (x) \xi^k \otimes \xi^k  \, .
\end{equation}
We can proceed by adding successively the $n(n+1)/2$ spiraling perturbations given in \eqref{e:spiral} corresponding to $\xi^k$ and amplitude $a_{q,k}(x)$. Observe that from \eqref{e:decomp_local} we have $\|a_{q,k}\|_0\leq \|h_q\|_0^{\sfrac12}\leq 2\delta_{q+1}^{\sfrac12}$. Since 
\[
Du_q^TDu_q+h_q=g_{q+1},
\]
we obtain a new immersion $u_{q+1}:U\to \R^N$ such that
\begin{equation}\label{e:iterate1}
\begin{split}	
\|u_{q+1}-u_q\|_0&\leq \delta_{q+1}\,\\
\|g_{q+1}-Du_{q+1}^TDu_{q+1}\|_0&\leq c_0\delta_{q+2}\,\\
[v_{q+1} - v_q]_1 &\leq C\delta_{q+1}^{\sfrac12}\, .
\end{split}
\end{equation}
From these estimates we easily conclude the $C^1$ convergence to an isometry.

\subsection{The quadratic term wins}\label{s:quadratic}

From a PDE point of view it is interesting to take a closer look at the calculation leading to \eqref{e:newmetric}, which forms the basis of the iteration above. Let us write \eqref{e:spiral} as 
\[
\tilde u(x)=u(x)+w(x),
\]
so that the new metric has the form:
\begin{equation}\label{e:expansion}
\partial_i\tilde u\cdot\partial_j\tilde u  = \partial_iu\cdot\partial_ju+\underbrace{(\partial_iu\cdot\partial_jw+\partial_ju\cdot\partial_iw)}_{=: L} + \underbrace{\partial_iw\cdot\partial_jw}_{=:Q}\, .
\end{equation}
The decomposition above simply gives the perturbation induced in the metric tensor by the perturbing map $w$ as a sum of the parts which
are, respectively, linear and quadratic in $w$.  
Recalling the orthogonality conditions \eqref{e:normalvectors} we see that $w$ is orthogonal to $\partial_i u$ for all $i$. Therefore
\[
\partial_iu\cdot\partial_jw =\partial_j(\partial_iu\cdot w)-\partial_i\partial_ju\cdot w=-\partial_i\partial_ju\cdot w,
\] 
so that
\begin{equation}\label{e:linear_estimate}
\|L\|_0 \leq C[u]_2\|w\|_0=O(\lambda^{-1})\, .
\end{equation}
On the other hand 
\begin{equation}\label{e:quadratic}
Q= a^2 \left(\cos^2( \lambda \xi\cdot x) + \sin^2 (\lambda \xi\cdot x)\right)\, \xi_i\xi_j +O(\lambda^{-1})= a^2 \xi_i\xi_j+O(\lambda^{-1})\, .
\end{equation}
We see that the specific oscillatory form of the perturbation $w$ in \eqref{e:spiral} makes the quadratic part much more important than the linear one: this seems a rather ``odd'' approach from a classical PDE point of view. From \eqref{e:linear_estimate} we also see that along the iteration the underlying frequencies $\lambda_q$ need to converge to $+\infty$ very fast. In particular it is clear that along the iteration the second derivatives of the immersions diverge.

\section{$C^{1,\alpha}$ isometric maps}\label{s:Borisov}

As we have seen, the construction above cannot possibly produce isometric immersions which are $C^2$. In the specific case of the Weyl problem, where $n=2$ and $N=3$, this is of course not surprising in light of the classical rigidity results of Herglotz and Cohn-Vossen. An interesting question is to understand if and where there is a sharp border on the H\"older scale $C^{1,\theta}$, $\theta\in (0,1)$ between the dramatically different behavior of solutions of the Weyl problem for low versus high $\theta$. 

 In a series of papers in the 1950s, cf. \cite{Borisov58-1,Borisov58-2,Borisov58-3,Borisov58-4}, Yu.~Borisov showed that the rigidity of the Weyl problem can in fact be extended to $C^{1,\theta}$ immersions provided $\theta$ is sufficiently large. 

\begin{theorem}\label{t:2/3}
Let $(\mathbb S^2,g)$ be a surface with $C^2$ metric and positive Gauss curvature, and let
$u\in C^{1,\theta}(\mathbb S^2;\R^3)$ be an isometric immersion 
with $\theta>2/3$. Then $u (\mathbb S^2)$ is the boundary of an open convex set.
\end{theorem}

Borisov's Theorem is more general, but the statement above avoids the introduction of Pogorelov's concept of bounded extrinsic curvature, cf. \cite{ConDS}: Borisov proves such property without any assumption on the topology of the surface and then exploits the work of Pogorelov, \cite{Pogorelov73}, to conclude the local convexity of the image. We will discuss later (cf. Section \ref{s:speculations}) a more recent, very short, proof of Borisov's Theorem discovered in \cite{ConDS}, which exploits the same key computation of Constantin-E-Titi's proof of part (a) of Onsager's conjecture: another remarkable analogy with Problem \ref{p:2d-Hoelder}!

On the other hand for sufficiently small H\"older exponents the Nash-Kuiper construction remains valid:
\begin{theorem}\label{t:1/7}
Let $(\Sigma, g)$ be a $C^2$ Riemannian manifold of dimension $n$. Any short immersion $u: \Sigma \to \mathbb R^{n+1}$ can be uniformly approximated with $C^{1,\theta}$ isometric immersions with
\begin{itemize} 
\item[(a)] $\theta < \frac{1}{1+n(n+1)}$ when $\Sigma$ is a closed ball;
\item[(b)] $\theta < \frac{1}{1+ n (n+1)^2}$ when $\Sigma$ is a general compact $n$-manifold.
\end{itemize}
The maps can be chosen to be embeddings if $u$ is an embedding. 
\end{theorem}

Case (a) of this theorem was announced in \cite{Borisov1965} by Yu.~Borisov, based on his habilitation thesis, under the additional assumption that $g$ be analytic.  A proof with $n=2$ appeared more than 40 years later, cf. \cite{Borisov2004}. The general statement of Theorem \ref{t:1/7} has been proved in \cite{ConDS}. 

We will discuss below the most relevant aspects of the argument and, in particular, the significance of the thresholds in (a) and (b). Observe that in the first interesting case of $2$-dimensional disks we have $\frac{1}{7}$: there is thus a significant gap between this and the ``rigidity threshold'' $\frac{2}{3}$ in Theorem \ref{t:2/3}. It is of course very tempting to ask whether there is a single sharp interface distinguishing between the two behaviors. Gromov in his article mentions $\frac{1}{2}$ (cf. Question 36 therein) as a possible threshold and we will discuss in Section \ref{s:speculations} some facts in favor of the latter conjecture. In the case of $2$-dimensional disks the very recent paper \cite{DIS} gave the first improvement of Borisov's local exponent, namely we have the following

\begin{theorem}\label{t:1/5}
Let $\overline{D}\subset \mathbb R^2$ be a closed disk and $g$ a $C^2$ metric on it. Then any short immersion $u: \overline{D}\to \mathbb R^3$ can be uniformly approximated with $C^{1,\theta}$ isometric immersions if $\theta < \frac{1}{5}$. The maps can be chosen to be embeddings if $u$ is an embedding. 
\end{theorem}

\subsection{The H\"older Nash iteration} 

Let $U\subset \R^n$ be an open domain with a smooth Riemannian metric $g$. Assume for the moment that we can carry the iteration as in Section \ref{s:Nashiteration} and consider once more the sequence of smooth immersions $u_q:(U,g_q)\to\R^N$ from Section \ref{s:Nashiteration}. Recall that
\begin{equation}\label{e:metric_error}
\|g_q-Du_q^TDu_q\|_0\leq c_0\delta_{q+1}
\end{equation}
where $g_q=g-\delta_q\Id$. The map $u_{q+1}$ is obtained by adding $s_n:=n(n+1)/2$ spiraling perturbations, so that
\[
u_{q+1} = u_q + \sum_{k=1}^{s_n} w_{q+1,k}\, ,
\]
each of the form
\[
w_{q+1,k}(x)=\frac{a_{q,k}(x)}{\lambda_{q+1,k}}\biggl(\sin(\lambda_{q+1,k} x\cdot\xi^k)\zeta^{q,k}(x)+\cos(\lambda_{q+1,k} x\cdot \xi^k)\eta^{q,k}(x)\biggr)\,,
\]
where the amplitudes $a_{q,k}$ are given by \eqref{e:decomp_local}, the unit vectors $\zeta^{q,k},\eta^{q,k}$ are normal to $u_{q,k-1}(U)$ and the frequencies $\lambda_{q,1}\leq \dots \leq \lambda_{q,s_n}$ still need to be chosen appropriately. 
For convenience set $\lambda_{q+1} = \max_k \lambda_{q+1,k} = \lambda_{q+1, s_n}$. 

As we have seen, 
\begin{equation}\label{e:amplitude}
\|a_{q,k}\|_0\leq \|g_q-Du_q^TDu_q\|_0^{\sfrac12}\leq \delta_{q+1}^{\sfrac12}.
\end{equation} 
Hence, neglecting lower order terms, we obtain
\begin{equation}\label{e:norme_Ck}
[u_{q+1} - u_q]_{1+m} \lesssim \delta_{q+1}^{\sfrac{1}{2}} \lambda_{q+1}^m \qquad \mbox{for $m\in \N$.}
\footnote{Here and in what follows, $A\lesssim B$ means that $A\leq cB$ for some universal constant $c$, and $A\lesssim B$ if $A\lesssim B$ and $B\lesssim A$.}
\end{equation}
By classical interpolation we conclude
\begin{equation}\label{e:C0alfa}
[Du_{q+1} - Du_q]_{\theta} \lesssim \delta_{q+1}^{\sfrac{1}{2}} \lambda_{q+1}^\theta\, ,
\end{equation}
where $[f]_{\theta}$ denotes the usual H\"older seminorm 
\[
\sup_{x,y\in U, x\neq y}\frac{|f(x)-f(y)|}{|x-y|^{\theta}}.
\]
The convergence in $C^{1,\theta}$ depends then on whether the sum $\sum_{q}\delta_{q}^{\sfrac{1}{2}} \lambda_{q}^\theta$ converges. 
In particular, if we can choose $\{\delta_q,\lambda_q\}_{q\in\N}$ such that 
\begin{equation}\label{e:Tristan}
	\lambda_q:=\lambda^{q}\quad\textrm{ and } \delta_q:=\lambda_q^{-2\theta_0}
\end{equation}
for some $\lambda>1$ and $\theta_0\in (0,1)$, then $\theta_0$ will be the threshold H\"older exponent for the convergence of the scheme.
We are confronted with two issues: we wish to have a fast convergence of $\delta_q$ to $0$ and a tame blow-up of $\lambda_q$. On the other hand the latter must be chosen large in order to make some errors negligible. 

To get an idea of whether such a choice of $\delta_q,\lambda_q$ is possible, recall the computation in Section \ref{s:quadratic} and in particular the estimate in \eqref{e:linear_estimate} for the error. Based on \eqref{e:amplitude} and \eqref{e:norme_Ck}\footnote{Since $\|u_{q+1}-u_q\|_2$ should blow-up, $\delta_q^{\sfrac{1}{2}} \lambda_q$ should also blow up and the exponential ansatz gives the estimate $\|u_q\|_2 \lesssim \delta_q^{\sfrac{1}{2}} \lambda_q$.} we expect\footnote{This is only the part of the metric error coming from the linear part $L$ of the perturbation. It can be checked that the error coming from the quadratic part $Q$ is smaller.} for $u_{q,1}:=u_q+w_{q+1}$
\begin{equation*}
\begin{split}
\|g_q+a_{q,1}^2\xi^{1}\otimes\xi^{1}-Du_{q,1}^TDu_{q,1}\|_0&\lesssim \frac{\|a_{q,1}\|_0[u_{q}]_2}{\lambda_{q+1,1}}\\
&\lesssim \frac{\delta_{q+1}\delta_{q}^{\sfrac12}\lambda_q}{\lambda_{q+1,1}}.
\end{split}	
\end{equation*}
At the second step we will however bring into play the second derivative of $u_{q,1}$, which we expect to be of size $\delta_{q+1}^{\sfrac{1}{2}} \lambda_{q+1,1}$. After $s_n$ steps, we can guess for $u_{q+1}=u_{q+1,k}$ an estimate of type
\begin{equation}\label{e:guess}
\|g_{q+1}-Du_{q+1}^TDu_{q+1}\|_0 \lesssim 
\delta_{q+1}^{\sfrac{1}{2}} \delta_q^{\sfrac{1}{2}} \lambda_q\lambda_{q+1,1}^{-1} + \sum_{k=2}^{s_n} \delta_{q+1} \lambda_{q+1,k-1} \lambda_{q,k}^{-1}\, . 
\end{equation}
This turns out to be correct, although the discussion above is somewhat simplified: there are several other error terms which must be computed.
Taking \eqref{e:guess} for granted, in order to keep \eqref{e:metric_error} we need to have
\begin{align*}
\delta_{q+2} \lesssim & \delta_{q+1}^{\sfrac{1}{2}} \delta_q^{\sfrac{1}{2}} \lambda_q \lambda_{q+1,1}^{-1}\\
\delta_{q+2} \lesssim & \delta_{q+1} \lambda_{q+1,k-1} \lambda_{q+1,k}^{-1}\,.
\end{align*}
If we optimize upon our choice of the parameters, these relations lead to 
\begin{equation}\label{e:Tristan2}
\delta_{q+2}^{s_n} \sim \delta_{q+1}^{s_n -1/2} \delta_q^{\sfrac{1}{2}} \lambda_q \lambda_{q+1}^{-1}\, .
\end{equation}
In view of \eqref{e:Tristan}, the latter identity takes the form
\[
\lambda^{(q+1) -2 s_n (q+2) \theta_0} \sim \lambda^{- (2s_n-1) (q+1) \theta_0 - \theta_0 q +q}\, .
\]
Taking the logarithm we easily conclude
\begin{equation}\label{e:exponent_iso}
\theta_0 = \frac{1}{1+2s_n}\, .
\end{equation}

\subsection{Borisov's exponents and beyond}
The somewhat exotic exponents of Theorem \ref{t:1/7} can now be easily explained: the threshold $\theta_0$ is related to the number of steps in a stage, i.e.~the number of spirals needed to obtain a full-rank correction of the metric error. In the case of Theorem \ref{t:1/7}(a) we can use a finite number of steps in the general Nash-Kuiper scheme to get a new short map from which we can proceed as in Section \ref{s:Nashiteration}: taking advantage of the minimal decomposition \eqref{e:linear_comb} we reach the threshold $\frac{1}{1+2s_n}= \frac{1}{1+n (n+1)}$.  In the case of Theorem \ref{t:1/7} (b) we have to use an additional partition of unity in $\Sigma$, and to control the overlaps of different charts requires a factor of $(n+1)$ more spirals. In other words each stage of the iteration consists of $(n+1)s_n=\frac{n(n+1)^2}{2}$ steps. Note that the general decomposition in \eqref{e:decomposition} would lead to even more steps and hence to a lower H\"older exponent.

 The heart of the matter in Theorem \ref{t:1/5} is that in $2$ dimensions we can hope to use, at each stage, a conformal transformation of coordinates that brings the metric error in diagonal form, thus allowing us to decompose it as the sum of $2$ rank-one terms rather than $3$. In the next section we will give a glimpse of some important technical obstructions for Theorem \ref{t:1/7}. The same obstructions appear in the proof of Theorem \ref{t:1/5} but the additional source of nontrivial complications is that the regularity of the conformal change of coordinates needed at each stage deteriorates dramatically as $q$ increases: it is apriori not even clear that such a scheme would at all converge in $C^1$. Indeed the related estimates do not allow to impose an exponential growth of the frequencies and we have to resort to a double exponential ansatz. 
 
It is not difficult to see that, if we enlarge the codimension, the argument of Theorem \ref{t:1/7} gives higher thresholds since then we can add spirals in parallel. In fact for $N =  n + s_n$ and $\Sigma$ equal to a  ball, a straightforward adaptation of the proof in \cite{ConDS} gives the threshold $\frac{1}{3}$. More work is needed to reach the threshold $\frac{1}{2}$ when $N$ is even higher, but this does not require any new insight. Instead a substantial new idea is needed to overcome $\frac{1}{2}$: this was achieved by A.~K\"allen in \cite{Kallen}. K\"allen's approach has the remarkable outcome that the Nash 1954 iteration scheme can be pushed ``almost'' up to $C^2$: for any metric of class $g\in C^{1,\theta}$ with any $\theta<1$ one obtains isometric embeddings of class $C^{1,\theta}$. Although very interesting from the PDE point of view, this result has less geometric impact:  the codimension needed is so large that it exceeds the one needed by Gromov to prove the existence of smooth isometric approximations when $g$ is smoother (we note in passing that it is still not known whether the $C^2$ regularity for $g$ is enough to show the existence of $C^2$ embeddings: the best result in that direction, due to Jacobowitz, needs a $C^{2,\beta}$ metric with positive $\beta$, cf. \cite{Jacobowitz}).

\subsection{Mollification and commutator estimate} The scheme outlined in the previous section has one drawback: there is a ``loss of derivative'' in the estimates. In particular observe that the perturbation $w_{q,1}$ involves taking vector fields normal to $u_q$ and thus depends certainly on the first derivative of $u_q$. This loss of derivatives propagates along the steps and stages of the iteration: the $j$-th derivative of $u_{q+1}$ depends certainly on the $j+1$-th, $j+2-th, \ldots$ and $j+s_n^*$-th derivatives of $u_q$, where $s^*n$ is the number of steps needed. This of course brings in higher and higher derivatives of the metric as well: it is in order to overcome this issue that Borisov assumes real analyticity of the metric $g$ (note that he also needs real analyticity for the starting short map $u_0$, but this can be assumed without loss of generality by a first regularizing procedure). 

As it is well-known in the PDE literature (following the other landmark work of Nash \cite{Nash1956}!), one way to overcome a loss of derivative in an iteration scheme is to introduce a mollification at each stage -- one may hope that this works provided the convergence rate is very fast. Thus, rather than defining $u_{q+1}$ in terms of $u_q$ we can define it in terms of $u_q* \varphi_\ell$, where $\varphi$ is a standard mollifier and $\ell = \ell_q$ a suitable mollification scale. The introduction of the latter scale is a real advantage only if $\ell_q \geq \lambda_q^{-1}$. However, under this assumption we have to ensure that $(u_q*\varphi_\ell)^\sharp e$ is close enough to $u_q^\sharp e$. In order to do this we exploit crucially two facts: the smallness of $g - u_q^\sharp e$ and an elementary commutator estimate between products and convolutions, which we state here.

\begin{lemma}\label{l:mollify}
Let  $\varphi\in C^\infty_c(\R^n)$ be symmetric and such that $\int\varphi=1$. Then for any $r\geq 0$ any $\theta\in (0,1]$ there
is a constant $C (r,\theta, n)$ such that the following estimate holds for any pair of $C^\theta$ functions $f$ and $g$: 
\begin{equation}\label{e:mollify3}
\|(fg)*\varphi_\ell-(f*\varphi_\ell)(g*\varphi_\ell)\|_r\leq C\ell^{2\theta -r}\|f\|_\theta\|g\|_\theta.
\end{equation}
\end{lemma}

This simple estimate due to Constantin, E and Titi has a very elementary proof and plays a crucial role both in their proof, cf. \cite{CET}, of the ``rigidity part'' of Onsager's Conjecture and in the short proof given in \cite{ConDS} of Borisov's Rigidity Theorem, cf. Section \ref{s:speculations} below.

\section{The Euler equations and Onsager's conjecture}\label{s:Euler}

The incompressible Euler equations describe the motion of a perfect incompressible fluid. Written down by L.~Euler over 250 years ago, these are the continuum equations corresponding to the conservation of momentum and mass of arbitrary fluid regions. In Eulerian variables they can be written as
\begin{equation}\label{e:Euler}
\left\{\begin{array}{l}
\partial_t v+(v\cdot \nabla) v+\nabla p =0\\ \\
\div v =0,
\end{array}\right.
\end{equation}
where $v=v(x,t)$ is the velocity and $p=p(x,t)$ is the pressure. We will focus on the 3-dimensional case 
with periodic boundary conditions. In other words we take the spatial domain to be the flat 3-dimensional torus $\T^3=\R^3/(2\pi\Z)^3$.

A classical solution on a given time interval $[0,T]$ is defined to be a pair $(v,p)\in C^1(\T^3\times [0,T])$. Despite the rich geometric structure underlying these equations (see e.g.~\cite{Constantin:2007tf} and references therein), little is known about classical solutions except
\begin{itemize}
\item[(i)] local well-posedness (i.e.~existence and uniqueness for short time) in H\"older spaces $C^{1,\theta}$, $\theta>0$ \cite{Lichtenstein:1968vc} or Sobolev spaces $H^s$, $s>5/2$ \cite{Ebin:1970vu,Kato:1972uo};
\item[(ii)] the celebrated blow-up criterion of Beale-Kato-Majda \cite{Beale:1984vk} and its geometrically refined variants, see e.g.~\cite{Constantin:1996bx}. 
\end{itemize}

\subsection{The paradox of Scheffer}
There are various notions of weak solutions (see for instance the survey article \cite{DS-Survey} and \cite{Sz-Roscoffnotes}), and despite the fact that uniqueness in general fails for such notions (see Theorem \ref{t:euler1} below and \cite{DS-ARMA,Daneri,DaneriSz} for further results), weak solutions have been studied because of their possible relevance to homogeneous 3D turbulence \cite{Onsager,Eyink,CET,CheskidovShvydkoy}. In particular we will consider pairs $(v,p):\T^3\times[0,1]\to\R^3\times\R$ which form a solution of \eqref{e:Euler} in the sense of distributions\footnote{Recall the classical computation that $(v\cdot \nabla )v=\div(v\otimes v)$ if $\div v=0$, so that distributional solutions are defined for any $v\in L^2(\T^3\times[0,1])$.}.

In contrast with the local well-posedness for classical solutions of \eqref{e:Euler}, weak solutions are in general quite ``wild'', and exhibit a behavior which is very different from classical solutions. Here we merely state Scheffer's amazing result from 1993 and refer to previous surveys \cite{DS-Survey,Sz-Roscoffnotes} for further results on distributional solutions.
\begin{theorem}\label{t:euler1}
There exist infinitely many non-trivial weak solutions $v\in L^{\infty}(\T^3\times\R)$ of \eqref{e:Euler} 
which have compact support in time.
\end{theorem}
This theorem was proved first by V.~Scheffer \cite{Scheffer93} in two dimensions for $v\in L^2 (\R^2\times\R)$. A.~Shnirelman \cite{Shnirelman97} subsequently gave a different proof for $v\in L^2(\T^2\times\R)$. The statement for arbitrary dimensions $d\geq 2$ and bounded velocities was obtained in $\mathbb R^d$ by \cite{DS-Annals}. 

\subsection{Energy conservation and Onsager's conjecture}

For classical solutions (i.e.~if $v\in C^1$) the total energy 
\[
E(t):=\frac{1}{2}\int_{\T^3}|v(x,t)|^2\,dx
\]
is conserved by the flow induced by \eqref{e:Euler}, so that $E(t)=E(0)$. However, for weak solutions this may not be true.  
Indeed, one of the cornerstones of three-dimensional turbulence is anomalous dissipation: it is an experimentally observed fact that the rate of energy dissipation in the vanishing viscosity limit (more precisely the infinite Reynolds number limit) stays above a certain non-zero constant.   This phenomenon is expected to arise from a mechanism of transporting energy from large to small scales, thereby leading to a cascade of energy. 

Assuming that a turbulent fluid is represented by a solution of the incompressible Navier-Stokes equations, in the vanishing viscosity limit one obtains the system \eqref{e:Euler}. Since classical solutions conserve the energy, in this (vaguely defined) limiting process one expects to find weak solutions of the Euler equations. It was L.~Onsager in 1949 \cite{Onsager} who first formulated the corresponding mathematical problem: is there a threshold between $C^0$ and $C^1$ regularity for energy conservation? Based on calculations in Fourier space, he formulated the statement in Conjecture \ref{c:Onsager} (in fact he had a non-rigorous proof of part (a)). 

Part (a) of the conjecture is fully resolved \cite{Eyink,CET}, whereas concerning part (b) substantial progresses have been made in the last five years, starting from \cite{DS-Inv}, although the full conjecture with threshold exponent $1/3$ remains an outstanding open problem. Having fixed a certain specific space of (at least $L^2$) functions $X$, these results can be classified in the following two categories:
\begin{enumerate}
\item[(A)] There exists a nontrivial weak solution $v\in X$ of \eqref{e:Euler} with compact support in time.
\item[(B)] Given any smooth positive function $E= E(t)>0$, there exists a weak solution $v\in X$ of  \eqref{e:Euler} with 
\begin{equation}\label{e:energy_id}
\frac{1}{2}\int |v(x,t)|^2\,dx=E (t)\quad\forall\,t.
\end{equation}
\end{enumerate}
Obviously both types lead to non-conservation of energy and would therefore conclude part (b) of Onsager's conjecture if proved for the space $X=L^{\infty}(0,1;C^{1/3-\eps}(\mathbb{T}^3))$. So far the best results are as follows.
\begin{theorem}\label{t:onsager}\hfill

\begin{itemize}
\item[(i)] Statement (A) is true for $X=L^1(0,1;C^{1/3-\eps}(\mathbb{T}^3))$.\footnote{$v\in L^1(0,1;C^{1/3-\eps}(\mathbb{T}^3))$ if and only if there exists an integrable function $A:(0,1) \to \mathbb R^+$ such that $|v (x,t) - v(x', t)|\leq A(t) |x-x'|^{\frac{1}{3}-\eps}$ for all $t,x,x'$.}
\item[(ii)] Statement (B) is true for $X=L^{\infty}(0,1;C^{1/5-\eps}(\mathbb{T}^3))$.
\end{itemize}
\end{theorem} 

Statement (B) has been shown for $X=L^\infty (0,1; C^{1/10-\eps})$ in \cite{DS-JEMS}, whereas P. Isett in his PhD thesis \cite{Isett} was the first to prove Statement (A) for $X= L^\infty (0,1; C^{1/5-\eps})$, thereby reaching the current best ``uniform'' H\"older exponent for Part (b) of Onsager's conjecture. Subsequently, T. Buckmaster, the two authors and P. Isett proved
Statement (B) for $X=L^\infty (0,1; C^{1/5-\eps})$ in \cite{BDIS}. Finally, Statement (A) for $X=L^1(0,1;C^{1/3-\eps}(\mathbb{T}^3))$
   has been proved recently in \cite{BDS}, based on a clever modification of the scheme by T. Buckmaster \cite{Buckmaster}. The basic construction underlying all these results was first introduced in \cite{DS-Inv}, and draws heavily on the 1954 paper of Nash \cite{Nash1954}. 

\subsection{Relaxation, subsolutions and h-principle}

Before explaining the basic construction of weak solutions for Theorem \ref{t:onsager} above, it is useful to look at some of the similarities between the systems \eqref{e:equations} and \eqref{e:Euler}. These similarities are based on the observation that both systems can be written as a differential inclusion. More precisely, both systems fit into the framework introduced by L.~Tartar in the context of compensated compactness in the 1970s \cite{Tartar,DS-Survey}, which amounts to separation into a linear system of conservation laws and non-linear pointwise constitutive relations. For the system of isometries \eqref{e:equations}, if we introduce the unknown $Z=Du$, this amounts locally to:
\begin{equation}\label{e:isometricT}
\textrm{curl }Z=0,\quad Z^TZ=g\, .
\end{equation}
Similarly, he Euler equations \eqref{e:Euler} can be written as
\begin{equation}\label{e:EulerT}
\left.\begin{array}{l}\partial_tv+\div u+\nabla q=0\\ \\ \div v=0\end{array}\right\}, \quad v\otimes v-u=\frac{2}{3}\bar{e}\, \Id.
\end{equation}
Here the ``state variable'' is the triple $Z=(v,u,q)$ with $u$ being a traceless symmetric $3\times 3$ matrix-valued function and $\bar{e}=\frac{1}{2}|v|^2$ is the kinetic energy density. 

One of the questions studied by Tartar was the weak closure of systems of the form \eqref{e:isometricT} or \eqref{e:EulerT}: to understand the relaxation of the constitutive relations when one considers a sequence of approximate solutions $Z_j$ converging weakly. For the system \eqref{e:isometricT} this is equivalent to the question we looked at in Section \ref{s:relaxation1}: short maps correspond to solutions of 
\begin{equation}\label{e:shortT}
\textrm{curl }Z=0,\quad Z^TZ\leq g.
\end{equation}
The analogous relaxation for \eqref{e:EulerT} has been computed in \cite{DS-ARMA}. It is given by
\begin{equation}\label{e:subsolutions}
\left.\begin{array}{l}\partial_tv+\div u+\nabla q=0\\ \\ \div v=0\end{array}\right\}, \quad v\otimes v-u\leq \frac{2}{3}\bar{e}\Id.
\end{equation}
Solutions of this system are called ``subsolutions'' of the Euler equations, and in \cite{DS-Annals,DS-ARMA} it was shown that any subsolution can be approximated weakly* in $L^{\infty}$ by bounded (but highly discontinuous) weak solutions of \eqref{e:Euler}. The construction is based on a well-known path in the literature for differential inclusions \cite{Cellina,Bressan,DacMarc,MullSver}, in particular it exploits the Baire category theorem (although one can give a proof using the alternative ``Lipschitz convex integration'' developed in \cite{MullSver})- we refer to \cite{DS-Survey} for a detailed exposition of this work. 

Dealing with merely bounded (i.e. $v\in L^{\infty}$) weak solutions of the Euler equations \eqref{e:Euler}  is somewhat reminiscent of dealing with Lipschitz\footnote{Recall that Lipschitz maps are differentiable almost everywhere, hence we mean here Lipschitz maps satisfying \eqref{e:equations} almost everywhere.} solutions for \eqref{e:equations}. As pointed out by Gromov \cite{Gromov86}, such maps need not be isometric in the sense of \eqref{e:true_isometry} and may in fact collapse entire submanifolds to a single point. Nevertheless, even if  \eqref{e:true_isometry} holds, the existence of a large class of Lipschitz isometries is much less surprising than the Nash-Kuiper theorem, since we are allowed to ``fold'' our Riemannian manifold. In this way one can even impose that the target has the same dimension as
the manifold, cf. \cite{KirSpSz}. Gromov in his article refers to the corresponding scheme as ``broken convex integration''. 

The notion of subsolutions of the Euler equations is closely connected to the Reynolds equations in classical turbulence theories. 
Let $v$  be a (deterministic or random turbulent) weak solution of \eqref{e:Euler} and consider a certain averaging process leading to the decomposition
\[
v=\overline{v}+w
\]
where $\overline{v}$ is the ``average'' and $w$ is the ``fluctuation''. The Euler equations \eqref{e:Euler} for $v$ transform into
\begin{equation}\label{e:averageE}
 \left\{\begin{aligned}
       &\partial_t \bar{v}+\div (\bar{v}\otimes \bar{v})+\nabla \bar{p}=-\div \bar{R}\\ \\
&\div \bar{v}=0
        \end{aligned}\right.
\end{equation}
where 
\begin{equation}\label{e:Rstress}
\bar{R}=\overline{v\otimes v}-\overline{v}\otimes\overline{v}=\overline{w\otimes w}.
\end{equation}
Being an average of positive semidefinite tensors, it is easy to see that $\bar{R}$ is positive semidefinite. The system \eqref{e:averageE} is equivalent to \eqref{e:subsolutions}. Indeed, given a subsolution $(\bar{v},\bar{u},\bar{q})$ define
\[
\bar{R}=\frac{2}{3}\bar{e}\,\Id-\bar{v}\otimes \bar{v}+\bar{u},\quad \bar{p}=\bar{q}-\frac{2}{3}\bar{e}.
\]
Then $\bar{R}$ is positive semidefinite and $(\bar{v},\bar{p},\bar{R})$ is a solution of \eqref{e:averageE}. 
Conversely, any solution $(\bar{v},\bar{p},\bar{R})$ of \eqref{e:averageE} with $\bar{R}\geq 0$ defines a subsolution $(\bar{v},\bar{u},\bar{q})$ with energy density
\begin{equation}\label{e:generalizedlocalenergy}
\bar{e}=\frac{1}{2}\left(\tr \bar{R}+|\bar{v}|^2\right)=\frac{1}{2}\tr(\bar{R}+\bar{v}\otimes \bar{v})
\end{equation}
by setting 
\[
\bar{u}=\bar{R}-\frac{2}{3}\bar{e}\,\Id+\bar{v}\otimes \bar{v},\quad \bar{q}=\bar{p}+\frac{2}{3}\bar{e}.
\]
In light of this interpretation of $\bar{R}$, it is natural to define the \emph{generalized energy tensor} of a subsolution $(\bar{v},\bar{p},\bar{R})$ to be the time-dependent tensor 
\begin{equation}\label{e:tensor}
\int_{\T^3}(\bar{v}\otimes \bar{v}+\bar{R})\,dx,
\end{equation}
and the associated \emph{generalized total energy} to be by its trace (cf.~\eqref{e:generalizedlocalenergy}):
\[
\overline{E}(t)=\frac{1}{2}\int_{\T^3}|\bar{v}|^2+\tr \bar{R}\,dx.
\]
Observe that the system \eqref{e:averageE} is highly under-determined. 
An important problem in the theory of turbulence is to obtain further restrictions on the tensor $\bar{R}$ in the form of constitutive (closure) relations. Thus an interesting question is whether there are additional constraints in the specific case where $\bar{R}$ arises -- in analogy with \eqref{e:Rstress} -- as a weak limit 
\begin{equation}\label{e:Reynolds}
\bar{R}=(\mathrm{w-lim}_{k\to\infty}v_k\otimes v_k)-\bar{v}\otimes \bar{v},
\end{equation}
where $v_k\rightharpoonup \bar{v}$ is a sequence of H\"older continuous weak solutions. Indeed, weak convergence has long been considered as a useful tool to study ``deterministic turbulence'' \cite{Lax:1990eg}. It follows from \cite{DS-ARMA,SzWiedemann} that no such constraints exist for $L^\infty$ weak solutions. It was recently shown in \cite{DaneriSz} that no additional constraints exist also for H\"older-continuous solutions and that therefore \emph{any} positive definite tensor can arise as \eqref{e:Reynolds} from $C^{\sfrac15-\eps}$-weak solutions of Euler:

\begin{theorem}[h-principle]\label{t:hprinciple}
Let $(\bar{v},\bar{p},\bar{R})$ be a smooth solution of \eqref{e:averageE} on $\T^3\times[0,T]$ such that $\bar{R}(x,t)$ is positive definite for all $x,t$. Then there exists for any $\theta<1/5$ a sequence $(v_k,p_k)$ of weak solutions of \eqref{e:Euler} such that 
\[
|v_k(x,t)-v_k(x',t)|\leq C_k|x-x'|^{\theta}\quad\textrm{ for all $x,x'$}
\] 
holds, 
\[
v_k\overset{*}{\rightharpoonup} \bar{v}\quad\textrm{ and }\quad v_k\otimes v_k\overset{*}{\rightharpoonup} \bar{v}\otimes \bar{v}+\bar{R}\quad\textrm{ in }L^\infty
\] 
uniformly in time and furthermore for all $t\in [0,T]$
\[
\int_{\T^3}v_k\otimes v_k\,dx=\int_{\T^3}(\bar{v}\otimes \bar{v}+\bar{R})\,dx.
\]
\end{theorem}

Theorem \ref{t:hprinciple} says that any smooth subsolution of the Euler equations can be weakly* approximated by H\"older-continuous weak solutions with given energy tensor. Observe that the uniform approximation in the Nash-Kuiper statement in Theorems \ref{t:main_C1_1} and \ref{t:1/7}, i.e.~convergence $u_j\to u$ in $C^0$ for a sequence of $C^1$ solutions of \eqref{e:equations}, can be equivalently stated as the weak* convergence in $L^\infty$ of a sequence $Z_j$ of H\"older-continuous solutions of \eqref{e:isometricT}. Therefore Theorem \ref{t:hprinciple} can be seen as the analogue of the Nash-Kuiper theorem.

\section{The Nash-scheme for the Euler equations}\label{s:Onsager}

In this and in the next section we review the key ideas leading to the proofs of Theorem \ref{t:onsager}. Although the basic scheme follows the one introduced in \cite{DS-Inv} by the authors, the presentation here uses crucial ideas that were introduced subsequently in the PhD Theses of T.~Buckmaster and of P.~Isett.

The construction of continuous and H\"older-continuous solutions of \eqref{e:Euler} follows the basic strategy of Nash in the sense that at each step of the iteration we add a highly oscillatory correction as the spiral in \eqref{e:spiral}. Note that both \eqref{e:Euler} and the equation of isometries \eqref{e:equations} are quadratic -- the oscillatory perturbation is chosen in such a way as to minimize the linearization and making the quadratic part of leading order (cf.~Section \ref{s:quadratic}). Then, a finite-dimensional decomposition of the error (cf.~\eqref{e:decomp_local}) is used to control the quadratic part. There are, however, two important differences:
\begin{itemize}
\item The linearization of \eqref{e:equations} is controlled easily by using the extra codimension(s) in the proof of Nash (cf.~the choice of perturbation in \eqref{e:spiral} being orthogonal to the previous image). For Euler, the linearization of \eqref{e:Euler} leads to a transport equation, which is very difficult to control over long times and seem to require a kind of CFL condition. This issue is still the main stumbling block in the full resolution of Onsager's conjecture and it will be examined in detail in the next section.
\item The exponent $1/3$ of Onsager's conjecture requires a sufficiently good correction of the error at each single step, whereas in the Nash iteration several steps ($s_n$ steps) are required -- this leads to the threshold exponent $(1+2s_n)^{-1}$ in Theorem \ref{t:1/7}. Consequently one-dimensional oscillations, as used in the Nash-Kuiper scheme and, more generally, in convex integration, cannot be 
used\footnote{However a ``multistep iteration'' using one-dimensional oscillation is possible in the case of Euler as well, as it has been recently shown by Isett and Vicol in \cite{IsettVicol}. This allows the authors to implement the iteration for a general class of active scalar equations, albeit leading to suboptimal H\"older exponents.}  for part (b) of Conjecture \ref{c:Onsager}. Thus, instead of convex integration, one needs to use special families of stationary flows as the replacement of \eqref{e:spiral} (compare \eqref{e:decomp_local} with \eqref{e:beltramireynolds} below).
\end{itemize}

The goal is to construct a sequence of subsolutions $(v_q,p_q,R_q)$, i.e.~solutions of 
\begin{equation}\label{e:Reynolds}
\left\{\begin{array}{l}
\partial_t v_q+{\rm div}\, v_q\otimes v_q+ \nabla p_q = - {\rm div}\, R_q\\ \\
\div v_q =0\, 
\end{array}\right.
\end{equation}
and iteratively remove the error. As a first observation note that if one is only interested in measuring the ``distance'' of a smooth pair $(v_q,p_q)$ from being a solution of \eqref{e:Euler}, then only the traceless part of $R_q$ is relevant: we can write
\[
R_q=\rho_q\Id+\mathring{R}_q,
\]
where $\mathring{R}_q$ is a traceless $3\times 3$ symmetric matrix, since $\div(\rho_q\Id)=\nabla\rho_q$. Hence if $\mathring{R}_q=0$ then $v_q$ is a solution of the Euler equations (perhaps with a different pressure). Recall that we also aim in Theorem \ref{t:onsager} at satisfying in addition \eqref{e:energy_id}. A natural analogy of the metric error $g-Du_q^TDu_q$ and the Reynolds stress can then be obtained 
by choosing a sequence $E_q=E_q(t)$ with $E_q(t)\to E(t)$ and setting
\begin{align*}
& \rho_q(t):=\frac{1}{3(2\pi)^3}\left(E_{q+1}(t)-\frac{1}{2}\int_{\T^3}|v_q(x,t)|^2\,dx\right),\\ \\
&R_q(x,t):=\rho_q(t)\Id+\mathring{R}_q(x,t)
\end{align*}
(cf.~\eqref{e:metric_q} and \eqref{e:newmetricerror}).
Thus, our approximations will consist of smooth solutions $(v_q,p_q,R_q)$ of \eqref{e:Reynolds} such that $\tr R_q$ is a function of time only,
and we will use $\|R_q\|_0$ to measure the distance of the pair $(v_q, p_q)$ from being a solution of \eqref{e:Euler}-\eqref{e:energy_id}. 

An important difference between the Reynolds stress and the metric error is that the latter is uniquely determined
from the metric $g$ and the short map $u$, whereas the tensor $\mathring{R}$ is not at all uniquely defined from \eqref{e:Reynolds}. However it is possible to select a good
``elliptic operator'' which solves the equations ${\rm div}\, \mathring R = f$. The relevant technical lemma
is the following one. 

\begin{lemma}[The operator $\textrm{div}^{-1}$]\label{l:div-1}
There exists a homogeneous Fourier-multiplier operator of order $-1$, denoted 
\[
\div^{-1}:C^\infty (\T^3; \R^3)\to C^\infty(\T^3;\S_0^{3\times 3})
\]
such that, for any $f\in C^\infty (\T^3; \R^3)$ with average $\fint_{\T^3}f=0$ we have
\begin{itemize}
\item[(a)] $\div^{-1} f(x)$ is a symmetric trace-free matrix for each $x\in \T^3$;
\item[(b)] $\div \div^{-1} f = f$.
\end{itemize}
\end{lemma}

\subsection{The approximating sequence and its size}  In analogy with the Nash construction our aim is to build a sequence of triples $(v_q, p_q, \mathring{R}_q)$ solving \eqref{e:Reynolds} which converge uniformly to a triple $(v,p,0)$. Actually in what follows we will mostly focus on the velocity $v$. The sequence will be achieved iteratively by adding a suitable perturbation to $v_q$ and $p_q$. We thus set 
\[
w_q= v_q-v_{q-1}.
\]
As in the Nash-Kuiper iteration, the size of $w_q$ will be controlled with two parameters. The {\em amplitude} $\delta_q$ bounds the $C^0$ norm:
\begin{align}\label{e:wq_C0}
\|w_q\|_0 \lesssim\; & \delta_q^{\sfrac{1}{2}}\, .
\end{align}
Up to negligible errors the Fourier transform of the perturbation $w_q$ will be localized in a shell centered around
a given {\em frequency} $\lambda_q$. Hence
\begin{align}\label{e:wq_C1}
\|\nabla w_q\|_0 \lesssim\; \delta_q^{\sfrac{1}{2}} \lambda_q\, .
\end{align}
Along the iteration we will have $\delta_q\to 0$ and $\lambda_q\to \infty$
at a rate that is at least exponential. For the sake of definiteness and in analogy with \eqref{e:Tristan} we may think
\begin{equation}\label{e:Tristan1}
	\lambda_q:=\lambda^{q}\quad\textrm{ and } \delta_q:=\lambda_q^{-2\theta_0}
\end{equation}
for some $\lambda>1$
(although in the actual proofs a slightly super-exponential growth is required).
Thus, as already discussed in the case of the Nash-Kuiper iteration, the positive number $\theta_0$ is the threshold H\"older regularity which we are able to achieve through the iteration. 

As in the Nash-Kuiper iteration, the perturbation $w_{q+1}$ is added to ``balance'' the error $R_q$: following the discussion of the previous section we can expect that $R_q \sim w_{q+1}\otimes w_{q+1}$ and for this reason we assume
\begin{align}
\|\mathring{R}_q\|_0 \leq\; & c_0\delta_{q+1}\label{e:est_R_C0}\\
\|\nabla \mathring{R}_q\|_0 \lesssim\; & \delta_{q+1} \lambda_q\, 
\end{align} 
for some small dimensional constant $c_0$ (in analogy with \eqref{e:metricerror_q}).
It turns out that along the iteration the perturbation $p_q-p_{q-1}$ behaves quadratically\footnote{This will lead to H\"older continuity of the pressure with exponent $2\theta$. Such an improvement in the H\"older exponent can also be obtained directly from  Schauder estimates for the pressure from the equation $-\Delta p=\div\div v\otimes v$. We learned about this improved Schauder estimate from L.~Silvestre first, but the same observation was also made independently by P.~Isett in \cite{Isett1}.}  in the perturbation $w_q$ and thus 
\begin{align}
\|p_q\|_0 \lesssim\; & \delta_{q}\\
\|\nabla p_q\|_0 \lesssim\; & \delta_{q} \lambda_q\, .
\end{align} 
So far we have not made any assumption on the size of the time derivatives. A key remark of Isett in \cite{Isett} compared to \cite{DS-Inv,DS-JEMS} is that advective derivatives behave much better than simple time derivatives. For instance, since 
\[
\partial_t v_q + (v_q\cdot \nabla) v_q = - \nabla p_q - {\rm div}\, R_q\, ,
\]
we have
\begin{equation}
\|\partial_t v_q + (v_q\cdot \nabla) v_q\|_0 \lesssim \delta_q\lambda_q\, .
\end{equation}
Note that, instead, $\|\partial_t v_q\|\lesssim \delta_q^{\sfrac{1}{2}} \lambda_q$. We have thus ``gained'' an extra factor $\delta_q^{\sfrac{1}{2}}$. 
The most important idea of Isett is that this gain holds also for the advective derivative of the Reynolds stress:
\begin{align}
\|\partial_t \mathring{R}_q + (v_q\cdot \nabla) \mathring{R}_q\|_0 \lesssim \delta_{q+1} \delta_q^{\sfrac{1}{2}} \lambda_q\, .\label{e:advective}
\end{align}
Finally, the control on the energy will be assumed to be of the following nature: set
\[
E_q(t)=(1-\delta_{q+1})E(t),
\]
so that 
\begin{equation}\label{e:energy_iter}
|\rho_q (t)| = \left|\frac{1}{3(2\pi)^3} \left(E_{q+1}(t) - \frac{1}{2}\int |v_q|^2 (x,t)\, dx \right)\right| \leq \frac{1}{4}\delta_{q+1}\, .
\end{equation}
The analogue of the strict shortness in the case of isometries is given by the positive definiteness of $R_q=\rho_q\Id+\mathring{R}_q$. 
However, estimates \eqref{e:est_R_C0} together with \eqref{e:energy_iter} lead to a stronger condition: it is easy to verify that for any given $r>0$, with a sufficiently small choice of $c_0$ (depending only on $r$) these two estimates ensure that
\begin{equation}\label{e:Rq}
R_q(x,t)\in \mathcal{C}_{r}\quad\textrm{ for all $(x,t)$}\quad \mbox{ and }\quad
\|R_q\|_0\sim \delta_{q+1}\, ,
\end{equation}
(recall the definition \eqref{e:cone} of the cone $\mathcal{C}_r$). 

\subsection{The oscillatory ansatz}  In analogy with Nash's approach to Proposition \ref{p:iter_C1} our strategy is to make $w_q$ a highly oscillatory vector field. Guided by the role Nash spirals in \eqref{e:spiral}, let us consider 
\begin{equation}\label{e:ansatz_1}
w_o(x,t)=W\Bigl(v_q (x,t),  R_q (x,t),\lambda_{q+1} x,\lambda_{q+1} t\Bigr)\, ,
\end{equation}
where $W$ is a function which we are going to specify next\footnote{The pressure $p_{q+1}$ will be defined similarly as $p_{q+1}= p_q + P (v_q, R_q, \lambda_{q+1} x, \lambda_{q+1} t)$, but we will not enter into the details in our discussion, since its role is anyway secondary.}. 

First of all, the oscillatory nature of the perturbation requires us to impose that $W$ is periodic in the variable $\xi\in \T^3$.
Next, observe that $v_{q+1}$ must satisfy the divergence-free condition ${\rm div}\, v_{q+1} =0$ and $v+w_o$ is not likely to fulfill this: we need to add a suitable correction $w_c$ in order to satisfy it.
Indeed a stronger analogy with the isometric embedding problem would be to consider first a vector potential for $v_q$, namely to write $v_q$ as $\nabla \times z_q$ for some smooth $z_q$. Subsequently we would like to perturb $z_q$ to a new 
\[
z_{q+1} (x,t) = z_q (x,t) + \frac{1}{\lambda_{q+1}} Z (v (x,t), R (x,t), \lambda_{q+1} x, \lambda_{q+1} t)\, .
\] 
Computing $v_{q+1} := \nabla \times z_{q+1}$ we get
\[
v_{q+1} (x,t) = v_q (x,t) + \underbrace{(\nabla_\xi \times Z) (v (x,t), \tilde R (x,t), \lambda x, \lambda t)}_{(P)} + O \left( \frac{1}{\lambda}\right)\, .
\]
The term (P) would correspond to $w_o$ if we were able to find a vector potential $Z$ for $W$ which is {\em periodic in $\xi$}. This requires $\div_\xi W = 0$ and $\langle W \rangle = 0$, where we use the notation $\langle, \rangle$ to denote the average in the $\xi$ variable.

Similar considerations (see for instance \cite{Sz-Roscoffnotes}) lead to the following set of conditions that we would like to impose on $W$:
\begin{itemize}
\item $\xi\mapsto W(v,R,\xi,\tau)$ is $2\pi$-periodic with vanishing average, i.e.
\begin{equation}\label{e:H1}\tag{H1}
\langle W\rangle := \frac{1}{(2\pi)^3}\int_{\T^3}W(v,R,\xi,\tau)\,d\xi=0;
\end{equation}
\item The average stress is given by $R$, i.e. 
\begin{equation}\label{e:H2}\tag{H2}
\langle W\otimes W\rangle = R
\end{equation}
for all $R\in \mathcal{C}_r$;
\item The ``cell problem'' is satisfied:
\begin{equation}\label{e:H3}\tag{H3}
\left\{\begin{array}{l}
\partial_{\tau} W+ v\cdot\nabla_{\xi}W+\div_{\xi} (W\otimes W) + \nabla_\xi P =0\\ \\
\div_\xi W = 0\, ,
\end{array}\right.
\end{equation} 
where $P=P(v,R,\xi,\tau)$ is a suitable pressure;
\item $W$ is smooth in all its variables and satisfies the estimates
\begin{equation}\label{e:H4}\tag{H4}
|W|\lesssim |R|^{1/2},\,|\partial_vW|\lesssim |R|^{1/2},\,|\partial_R W|\lesssim |R|^{-1/2}.
\end{equation}
\end{itemize}
Observe that \eqref{e:H2} corresponds to \eqref{e:Rstress} and  \eqref{e:Reynolds}, 
\eqref{e:H3} arises from plugging the oscillatory ansatz \eqref{e:ansatz_1} into Euler, and \eqref{e:H4} are estimates consistent with \eqref{e:H2}.  

As a consequence of \eqref{e:H1}-\eqref{e:H2} we obtain
\[
\int_{\T^3}|v_{q+1}|^2\,dx\sim \int_{\T^3}|v_{q}|^2\,dx+\int_{\T^3}\langle |W|^2\rangle \,dx=\int_{\T^3}|v_{q}|^2\,dx+3(2\pi)^3\rho_{q}(t),
\]
so that \eqref{e:energy_iter} can be ensured inductively. The main issues are therefore
\begin{itemize}
\item to show that indeed it is possible to send $\delta_q$ to $0$ as $q\uparrow \infty$ (so that the scheme converges) 
\item and to obtain a relation between $\delta_q$ and $\lambda_q$ in the form of \eqref{e:Tristan1}.
\end{itemize}
We will see that, if we were able to find a ``profile $W$ satisfying (H1)-(H2)-(H3)-(H4), then the iteration proposed so far would lead to a proof of the Onsager's conjecture.

\subsection{$\frac{1}{3}$-scheme}\label{s:1/3} 
Assuming the existence of a such a profile $W$, the next stress tensor $\mathring{R}_{q+1}$ 
would then be defined through
\begin{eqnarray}
\mathring{R}_{q+1}&=\;-&\div^{-1}\Bigl[\partial_tv_{q+1} +\div(v_{q+1} \otimes v_{q+1})+\nabla p_{q+1}\Bigr]\notag\\
&=\;-&\underbrace{\div^{-1}\Bigl[\partial_t w_{q+1} +v_q \cdot\nabla w_{q+1}\Bigr]}_{=:\mathring{R}_{q+1}^{(1)}}\notag\\
&\phantom{=}\; -& \underbrace{\div^{-1}\Bigl[\div(w_{q+1}\otimes w_{q+1}- R_q)+\nabla (p_{q+1}-p_q)\Bigr]}_{=:\mathring{R}_{q+1}^{(2)}}\notag\\
&\phantom{=}\;-& \underbrace{\div^{-1}\Bigl[w_{q+1}\cdot \nabla v_q\Bigr]}_{=:\mathring{R}_{q+1}^{(3)}}
\end{eqnarray}
where $\div^{-1}$ is the operator of order $-1$ from Lemma \ref{l:div-1}. Since we are assuming that the size of the corrector $w_c$ is negligible compared to $w_o$, we will discuss the corresponding terms where $w_o$ replaces $w_q$. 

First expand $W (v, R, \xi, \tau)$ as a Fourier series in $\xi$. We then could compute
\begin{equation}\label{e:estR3}
\mathring{R}^{(3)}=\div^{-1}\Bigl[w_o\cdot \nabla v_q\Bigr]=\div^{-1}\sum_{k\in\Z^3, k\neq 0}c_k(x,t)e^{i\lambda_{q+1} k\cdot x}\, ,
\end{equation}
where the coefficients $c_k (x,t)$ vary much slower than the rapidly oscillating exponentials. When we apply the operator $\div^{-1}$ we can therefore treat the $c_k$ as constants and gain a factor $\frac{1}{\lambda_{q+1}}$ in the outcome: a typically ``stationary phase argument''. Note that it is crucial that $c_0$ vanishes: this is in fact the content of condition (H1).

Using \eqref{e:H4} we can estimate the size of each term $c_k$ as
\[
\|c_k\|_0\lesssim \|W\|_0\|\nabla v_q\|_0\lesssim \|R_q\|_0^{1/2}\|\nabla v_q\|_0.
\]
Applying \eqref{e:wq_C1} and \eqref{e:Rq} we arrive at 
\begin{equation}\label{e:Nash_term}
\|\mathring{R}_{q+1}^{(3)}\|_0 \lesssim \frac{\delta_{q+1}^{\sfrac{1}{2}} \delta_q^{\sfrac{1}{2}} \lambda_q}{ \lambda_{q+1}}\, .
\end{equation}
Coming to the two remaining terms observe that one needs to differentiate the perturbation $w_o$ in $x$ and $t$, where there is a distinction between ``slow'' and ``fast'' derivatives -- we refer to ``fast derivatives'' if the term involves a factor of $\lambda_{q+1}$. For instance
\[
\partial_t W=\underbrace{\partial_vW\partial_tv_q+\partial_RW\partial_tR_q}_{\textrm{slow}}+\underbrace{\lambda_{q+1}\partial_\tau W}_{\textrm{fast}}.
\]
Owing to condition \eqref{e:H3} (the ``cell problem'') the fast derivatives in $\mathring{R}_{q+1}^{(1)}+\mathring{R}_{q+1}^{(2)}$
vanish identically. Hence, by some abuse of notation, we may write
\begin{eqnarray}
\mathring{R}_{q+1}^{(1)}&=&\div^{-1}\Bigl[(\partial_t+v_q\cdot\nabla)^{\textrm{\tiny slow}}W\Bigr],\label{e:R-1}\\
\mathring{R}_{q+1}^{(2)}&=&\div^{-1}\Bigl[\div^{\textrm{\tiny slow}}(W\otimes W-R_q)\Bigr].\label{e:R-2}
\end{eqnarray}
Observe that the expression in \eqref{e:R-1} is linear in $W$, hence the same stationary phase argument as above applies. We calculate:
\[
(\partial_t+v_q\cdot\nabla)^{\textrm{\tiny slow}}W=\partial_vW\,(\partial_t+v_q\cdot\nabla)v_q+\partial_RW\, (\partial_t+v_q\cdot\nabla)R_q
\]
so that, writing as before,
\[
\mathring{R}_{q+1}^{(1)}=\div^{-1}\sum_{k\in\Z^3, k\neq 0}c'_k(x,t)e^{i\lambda_{q+1} k\cdot x}\,
\]
for some $c'_k$. This time, using \eqref{e:H4}, we have 
\[
\|c'_k\|_0\lesssim \|R_q\|_0^{1/2}\|(\partial_t+v_q\cdot\nabla)v_q\|_0+\|R_q\|_0^{-1/2}\|(\partial_t+v_q\cdot\nabla)R_q\|_0.
\]
From \eqref{e:wq_C1}, \eqref{e:advective} and  \eqref{e:Rq} we then deduce
\begin{equation*}
\|\mathring{R}_{q+1}^{(1)}\|_0 \lesssim \frac{1}{\lambda_{q+1}}\left(\delta_{q+1}^{1/2}\delta_q\lambda_q+\delta_{q+1}^{1/2}\delta_q^{1/2}\lambda_q\right)
\; \lesssim \frac{\delta_{q+1}^{1/2}\delta_q^{1/2}\lambda_q}{\lambda_{q+1}}\,.
\end{equation*}
Finally, observe that in \eqref{e:R-2} we have $\langle W\otimes W\rangle=R_q$ because of condition \eqref{e:H2}, so that once more, in the expansion
of $W\otimes W-R_q$ as a Fourier-series in $\xi$ there is no term $k=0$. Hence the same stationary phase estimate can be applied once more. Writing
\[
\mathring{R}_{q+1}^{(2)}=\div^{-1}\sum_{k\in\Z^3, k\neq 0}c''_k(x,t)e^{i\lambda_{q+1} k\cdot x}
\]
and using \eqref{e:H4} we have the estimate
\begin{equation*}
\begin{split}
\|c''_k\|_0&\lesssim \|W\|_0\|\partial_vW\|_0\|Dv_q\|_0+\|W\|_0\|\partial_RW\|_0\|D R_q\|_0\\
&\lesssim\|R_q\|_0\|D v_q\|_0+\|D R_q\|_0,
\end{split}
\end{equation*}
so that 
\begin{equation}\label{e:R2new}
\begin{split}
\|\mathring{R}_{q+1}^{(2)}\|_0&\lesssim \frac{1}{\lambda_{q+1}}\left(\delta_{q+1}\delta_q^{1/2}\lambda_q+\delta_{q+1}\lambda_q\right)\\
&\lesssim \frac{\delta_{q+1}\lambda_q}{\lambda_{q+1}}\,.
\end{split}
\end{equation}
Summarizing, we obtain 
\begin{equation}\label{e:Rnew}
\|\mathring{R}_{q+1}\|_0\lesssim \frac{\delta_{q+1}^{1/2}\delta_q^{1/2}\lambda_q}{\lambda_{q+1}}.
\end{equation}
Of course, this is just one of the estimates for $(v_{q+1},p_{q+1},R_{q+1})$ in \eqref{e:wq_C0}-\eqref{e:Rq}, similar ones should be obtained for all the other quantities. However, \eqref{e:Rnew} already implies a relation between $\delta_q$ and $\lambda_q$. Indeed, comparing it with \eqref{e:est_R_C0}, the inductive step requires 
\[
\delta_{q+2}\sim\frac{\delta_{q+1}^{1/2}\delta_q^{1/2}\lambda_q}{\lambda_{q+1}}.
\]
Assuming $\lambda_q\sim \lambda^q$ for some fixed $\lambda\gg 1$, this would lead to 
\begin{equation}\label{e:1/3}
\delta_q^{1/2}\sim \lambda^{-\sfrac{q}{3}}\sim \lambda_q^{-\sfrac{1}{3}},
\end{equation}
which, comparing with \eqref{e:Tristan2}, gives $\theta_0=1/3$ as the critical H\"older regularity. 

In the derivation above we have assumed the existence of $W$ with properties \eqref{e:H1}-\eqref{e:H4}. Next we will discuss how one could construct such $W$. As it turns out we are not able to fulfill all the conditions without further modifications. These modifications will eventually lead to additional error terms and are responsible for the lower threshold exponent $\theta_0=1/5$ in Theorem \ref{t:onsager}(ii).

\subsection{Beltrami flows} In this section we show how almost all conditions on the function $W=W(v,R,\xi,\tau)$ can be fulfilled. 
Let us first examine the simple case in which we set $v=0$: it is then possible to construct a function $W_s (R, \xi) = W (0, R, \xi, \tau)$ satisfying the constraints  \eqref{e:H1}-\eqref{e:H4}. The basic building block is given by Beltrami flows, which form the counterpart of the Nash spirals. Start with the identity
\[
\div(U\otimes U)=U\times\mathrm{curl}\,U-\tfrac12\nabla |U|^2\, ,
\]
for smooth $3$-dimensional vector fields $U$.
In particular any eigenspace of the curl operator (i.e.~the solution space of the system
\begin{equation*}
\left\{
\begin{array}{lll}
\mathrm{curl}\,U&=&\lambda_0 U\\ \\
\div U &=&0
\end{array}\right.
\end{equation*}
for $\lambda_0$ constant)
leads to a \emph{linear} space of stationary flows of the incompressible Euler equations. These can be written as
\begin{equation}\label{e:Bflow1}
\sum_{|k|=\lambda_0}a_kB_ke^{ik\cdot\xi}
\end{equation}
for normalized complex vectors $B_k\in\C^3$ satisfying 
\[
|B_k|=1,\quad k\cdot B_k=0\quad\textrm{ and }\quad ik\times B_k=\lambda_0 B_k,
\]
and arbitrary coefficients $a_k\in\C$. Choosing $B_{-k}=-\overline{B_k}$ and $a_{-k}=\overline{a_k}$ ensures that $U$ is real-valued. A calculation then shows
\begin{equation}\label{e:beltramireynolds}
\langle U\otimes U\rangle=\frac{1}{2}\sum_{|k|=\lambda_0}|a_k|^2\Bigl(\Id-\frac{k\otimes k}{|k|^2}\Bigr).
\end{equation}
Moreover, recalling the condition that $W$ must be $2\pi$-periodic in the $\xi$ variable, we impose that
$k\in \mathbb Z^3$.
The identity \eqref{e:beltramireynolds} leads to the following decomposition Lemma which is the analogue of Lemma \ref{l:decomp1}.

\begin{lemma}\label{l:Bdecomposition}
For every $N\in\N$ we can choose $0<r_0<1$ and $\lambda_0 > 1$ with the following property.
There exist pairwise disjoint subsets 
\[
\Lambda_j\subset\{k\in \Z^3:\,|k|=\lambda_0\} \qquad j\in \{1, \ldots, N\}
\]
and smooth positive functions 
\[
\gamma^{(j)}_k\in C^{\infty}\left(B_{r_0} (\Id)\right) \qquad j\in \{1,\dots, N\}, k\in\Lambda_j
\]
such that
\begin{itemize}
\item[(a)] $k\in \Lambda_j$ implies $-k\in \Lambda_j$ and $\gamma^{(j)}_k = \gamma^{(j)}_{-k}$;
\item[(b)] For each $R\in B_{r_0} (\Id)$ we have the identity
\begin{equation}\label{e:split}
R = \frac{1}{2} \sum_{k\in\Lambda_j} \left(\gamma^{(j)}_k(R)\right)^2 \left(\Id - \frac{k}{|k|}\otimes \frac{k}{|k|}\right) 
\qquad \forall R\in B_{r_0}(\Id)\, .
\end{equation}
\end{itemize}
\end{lemma}

This lemma, taken from \cite{DS-Inv} (see also \cite{Isett} for a geometric proof) allows us to choose the amplitudes as 
\begin{equation}\label{e:Bflow2}
a_k= \sqrt{{\rm tr}\, R}\,\gamma^{(j)}_k\left(\frac{R}{{\tfrac{1}{3}\tr\,} R}\right)
\end{equation}
for any $R\in \mathcal{C}_{r_0}$. 
With this choice of $a_k=a_k (R)$, we can then set 
\[
W_s (R, \xi):= \sum_{k\in \Lambda^{(1)}} a_k (R) B_k e^{i  k \cdot \xi}
\] 
(defined through the Beltrami-flow relation \eqref{e:Bflow1}). 
Note that that for such $W_s$ the sizes of $W$ and of any $R$-derivative of $W$ satisfy estimates \eqref{e:H4}.

\section{The transport problem and the $\frac{1}{5}$ threshold}\label{s:transport}

Having obtained a profile $W (0, R, \xi, \tau) = W_s (R, \xi)$, it seems natural to extend $W$ by imposing that $\partial_\tau W + v\cdot \nabla_\xi W =0$, leading to the formula
\begin{equation}
W (v, R, \xi, \tau) = W_s (R, \xi-v\tau)= \sum_{k\in \Lambda^{(1)}} a_k (R) B_k e^{i  (k - v \tau)\cdot}\, .
\end{equation}
However the latter fails to satisfy \eqref{e:H4}, because $|\partial_v W(v,R, \xi, \tau)|\sim |R|^{\sfrac12}|\tau|$. This is a serious problem: observing that $\tau$ is the ``fast time'' variable, in the construction \eqref{e:ansatz_1} $\tau=\lambda_{q+1}t$, leading to an additional factor $\lambda_{q+1}$ in the estimates for $\mathring{R}_{q+1}^{(1)}$ and $\mathring{R}_{q+1}^{(2)}$: this loss destroyes any hope that the scheme might converge.

In \cite{DS-Inv,DS-JEMS} a ``phase function'' $\phi_k(v,\tau)$ was introduced to deal with the transport part of the cell problem. By considering $W$ of the form
\begin{equation}\label{e:ansatz1}
\sum_{|k|=\lambda_0}a_k(R)\phi_k(v,\tau)B_ke^{ik\cdot\xi}
\end{equation}
the cell problem in \eqref{e:H3} leads to the equation
\[
\partial_\tau\phi_k+i(v\cdot k)\phi_k=0\, .
\]
Since the exact solution $\phi_k(v,\tau)=e^{-i(v\cdot k)\tau}$ is incompatible with the requirement \eqref{e:H4}, an approximation is used\footnote{To be precise, the approximation involves a partition of unity over the space of velocities and the use of $8$ distinct families $\Lambda^{(j)}$ in Lemma \ref{l:Bdecomposition}.} such that
\[
\partial_\tau\phi_k+i(v\cdot k)\phi_k=O\left(\mu_q^{-1}\right),\qquad |\partial_v\phi_k|\lesssim \mu_q
\]
for some new parameter $\mu_q$. This leads to the following corrections to \eqref{e:H3} and \eqref{e:H4}: \eqref{e:H3} is only satisfied approximately:
\begin{equation*}
\partial_{\tau} W+ v\cdot\nabla_{\xi}W+\div_{\xi} (W\otimes W) + \nabla_\xi P = O(\mu_q^{-1})
\end{equation*}
and in \eqref{e:H4} the second inequality is replaced by 
\begin{equation*}
|\partial_vW|\lesssim \mu_q|R|^{1/2}.
\end{equation*}

\subsection{Flow and CFL condition}
A further improvement was obtained in  \cite{BDIS}, following an idea first introduced by Isett in \cite{Isett}. We change the ansatz \eqref{e:ansatz1} on $W$ and look for a perturbation $w_o$ which has the form
\begin{equation}\label{e:Phil}
w_o (x,t) = W_s (R_q (x,t), \lambda_{q+1} \Phi_q (x,t)) = \sum_{k\in \Lambda^{(1)}} a_k (R_q (x,t)) B_k e^{i  \lambda_{q+1} \Phi_q (x,t)}\, , 
\end{equation}
where $\Phi_q$ solves the transport equation 
\begin{equation}\label{e:inverseflow}
\partial_t \Phi_q + (v_q \cdot \nabla_x) \Phi_q =0\, .
\end{equation}
With \eqref{e:Phil}, we would have
\begin{equation}\label{e:noncutoff-1}
\mathring{R}_{q+1}^{(1)} = \sum_{k\in \Lambda^{(1)}} \nabla a_k (R_q) (\partial_t R_q + (v_q\cdot \nabla) R_q) e^{i\lambda_{q+1} \Phi_q}
\end{equation}
and, assuming that $D\Phi_q (x,t)$ is not too far from the identity, the stationary phase argument together with the
bound \eqref{e:advective} would lead to 
\begin{equation}\label{e:non_cutoff}
\|\mathring{R}^{(1)}\|_0 \lesssim \delta_{q+1}^{\sfrac{3}{2}} \delta_q^{\sfrac{1}{2}} \lambda_q \lambda_{q+1}^{-1}\, .
\end{equation}
However, since $\|Dv_q\|_0\to \infty$, we expect the deformation matrix $D\Phi_q$ to be controllable only for short times. More precisely, by a well-known elementary estimate on ODEs, if $\Phi_q (x, t_0) =x$, then
\begin{equation}\label{e:CFL}
\|D\Phi_q (\cdot, t) - {\rm Id}\|_0 \lesssim \|\nabla v_q\|_0 |t-t_0| \lesssim \delta_q^{\sfrac{1}{2}} \lambda_q |t-t_0|
\end{equation}
for $|t-t_0|\lesssim (\delta_q^{\sfrac{1}{2}} \lambda_q)^{-1}$. The latter is a typical ``CFL condition'', cf. \cite{CFL}. 

To handle this problem we proceed as in \cite{BDIS} and consider a partition of unity $(\chi_j)_j$ on the time interval $[0,T]$ such that the support of each $\chi_j$ is an interval $I_j$ of size $\frac{1}{\mu_q}$ for some $\mu_q\gg 1$. In each time interval $I_j$ we set $\Phi_{q,j}$ to be the solution of the transport equation \eqref{e:inverseflow} which satisfies 
\[
\Phi_{q,j}(x,t_j)=x,
\] 
where $t_j$ is the center of the interval $I_j$. Recalling that $\|Dv_q\|_0\lesssim \delta_q^{\sfrac12}\lambda_q$, \eqref{e:CFL} leads to 
\begin{equation}\label{e:DPhi}
\|D\Phi_{q,j}\|_0=O(1)\quad\textrm{ and }\quad \|D\Phi_{q,j}-\Id\|_0\lesssim \frac{\delta_q^{\sfrac12}\lambda_q}{\mu_q}
\end{equation}
provided 
\begin{equation}\label{e:CFLcondition}
\mu_q\geq \delta_q^{\sfrac12}\lambda_q,
\end{equation}
an estimate we will henceforth assume. Observe also that $|\partial_t\chi_j|\lesssim \mu_q$.
 
The new fluctuation will take the form
\begin{align}
w_o &=\sum_j\chi_j(t)\sum_{k\in \Lambda^{(i(j))} }a_k (R_q)B_ke^{i\lambda_{q+1}k\cdot\Phi_{q,j}}\\
&=\sum_{k,j} a_{k,j} (R_q)\phi_{kj} B_ke^{i\lambda_{q+1}k\cdot x}\, ,
\end{align}
where:
\begin{itemize}
\item $i(j)$ equals $1$ if $j$ is odd and $2$ if $j$ is even;
\item $\Lambda^{(1)}$ and $\Lambda^{(2)}$ are two disjoint families of vectors from Lemma \ref{l:Bdecomposition};
\item the phase functions $\phi_{kj}$ are given by
$e^{i \lambda_{q+1} k \cdot (\Phi_j (x,t) -x)}$. 
\end{itemize}
In computing now $\mathring{R}_{q+1}^{(1)}$ we get, compared to \eqref{e:noncutoff-1}, an additional term of
the form 
\[
{\rm div}^{-1} \left[ \sum_j \partial_t \chi_j(t)\sum_{k\in \Lambda^{(i(j))} } a_k (R_q)\phi_{kj} B_ke^{i\lambda_{q+1}k\cdot x}\right]\, 
\]
and in view of $|\partial_t\chi_j|\lesssim \mu_q$ the estimate \eqref{e:non_cutoff} becomes
\begin{equation}\label{e:transport}
\|\mathring{R}_{q+1}^{(1)}\|_0 \lesssim \delta_{q+1}^{\sfrac{3}{2}} \delta_q^{\sfrac{1}{2}} \lambda_q \lambda_{q+1}^{-1} + \delta_{q+1}^{\sfrac{1}{2}} \mu_q \lambda_{q+1}^{-1}\stackrel{\eqref{e:CFLcondition}}{\lesssim} \delta_{q+1}^{\sfrac{1}{2}} \mu_q \lambda_{q+1}^{-1}\, .
\end{equation}
As for $\mathring{R}_{q+1}^{(3)}$ we can assume that \eqref{e:Nash_term} still holds. On the other hand the estimate for $\mathring{R}_{q+1}^{(2)}$ involves certainly some new error terms. First of all, since the profile $W_s$ solves
${\rm div}_\xi (W_s\otimes W_s) + \nabla_\xi P = 0$, there are no ``fast derivatives'' in the expression for $\mathring{R}_{q+1}^{(2)}$. Hence
\begin{equation}\label{e:overlapslow}
\mathring{R}^{(2)} = \div^{-1}\left[\div^{\textrm{\tiny slow}}(w_o\otimes w_o-R_q)\right]\, .
\end{equation}
We next compute
\begin{equation}\label{e:average-1/5}
\begin{split}
w_o\otimes w_o&= \frac{1}{2}\sum_{k,j} \chi_j^2|a_{k,j}|^2\left(\Id-\frac{k\otimes k}{|k|^2}\right)+\\&\quad+\sum_{j,j',k+k'\neq 0}\chi_j\chi_{j'}a_{kj}a_{k'j'}\phi_{kj}\phi_{k'j'}B_k\otimes B_{k'}e^{i\lambda_{q+1}(k+k')\cdot x}\\
&=R_q+\sum_{k''\neq 0}c_{k''}(x,t)e^{i\lambda_{q+1}k''\cdot x}\,.
\end{split}
\end{equation}
Since $\|\nabla \phi_{kj}\|_0 \leq \delta_q^{\sfrac{1}{2}} \lambda_q \lambda_{q+1} \mu_q^{-1}$ according to \eqref{e:DPhi}, we can estimate 
\[
\|\nabla c_{k''}\|_0 \lesssim \delta_{q+1} \lambda_q + \delta_{q+1} \delta_q^{\sfrac{1}{2}} \lambda_q \lambda_{q+1} \mu_q^{-1}\, .
\] 
Hence, by the stationary phase estimate, we expect
\begin{equation}\label{e:optimize}
\|\mathring{R}_{q+1}^{(2)}\|_0 \lesssim \delta_{q+1} \lambda_q \lambda_{q+1}^{-1} + \delta_{q+1} \delta_q^{\sfrac{1}{2}} \lambda_q \mu_q^{-1}\, .
\end{equation}
Combining \eqref{e:Nash_term},\eqref{e:transport} and \eqref{e:optimize} (and taking into consideration \eqref{e:CFLcondition}) we conclude
\begin{equation}\label{e:optimize2}
\|\mathring{R}_{q+1}\|_0 \lesssim \delta_{q+1}^{\sfrac{1}{2}} \mu_q \lambda_{q+1}^{-1}  + \delta_{q+1} \delta_q^{\sfrac{1}{2}} \lambda_q \mu_q^{-1}
\end{equation}
Optimizing in $\mu_q$ we then reach
\begin{equation}\label{e:optimized}
\|\mathring{R}_{q+1}\|_0 \lesssim \delta_{q+1}^{\sfrac{3}{4}} \delta_q^{\sfrac{1}{4}} \lambda_q^{\sfrac{1}{2}} \lambda_{q+1}^{-\sfrac{1}{2}}\, ,
\end{equation}
namely to
\[
\delta_{q+2} \sim \delta_{q+1}^{\sfrac{3}{4}} \delta_q^{\sfrac{1}{4}} \lambda_q^{\sfrac{1}{2}} \lambda_{q+1}^{-\sfrac{1}{2}}\, .
\]
Plugging \eqref{e:Tristan2} in the latter identity and taking logarithms leads to $\theta_0 = \frac{1}{5}$. 

\section{H-principle for H\"older solutions of Euler}\label{s:h-principle}

The Beltrami flows together with the transport ansatz explained in the previous subsections settle the issue of convergence (at least for H\"older exponents $\theta<1/5$), but are not sufficient to conclude the h-principle statement of Theorem \ref{t:hprinciple}. Indeed, the problem is reminiscent of the difference between the global form of the Nash stage in Proposition \ref{p:iter_C1} (which is based on the global decomposition \eqref{e:decomposition}) and the local version suitable for the iteration based on  Lemma \ref{l:decomp1}. It turns out that even when we increase the number of modes the Beltrami flows cannot generate arbitrary positive definite stresses (in other words the expression for $\langle U\otimes U\rangle$ in \eqref{e:beltramireynolds} cannot be an arbitrary positive definite matrix $R$; the set of possible $R$ which can be generated has been computed in \cite{Choffrut}). 

Nevertheless, there is a very simple set of stationary flows (which we will call ``Mikado flows'') based on pipe flow, which can generate all $R$. These flows were introduced in \cite{DaneriSz}. 
\begin{lemma}\label{l:Mikado}For any compact subset $\mathcal N$ consisting of positive definite $3\times 3$ matrices
there exists a smooth vector field 
\[
W:\mathcal N\times \T^3\to \R^3,\quad i=1,2
\]
such that, for every $R\in\mathcal N$ 
\begin{equation}\label{e:Mikado}
\left\{\begin{array}{l}
\div_\xi(W(R,\xi)\otimes W(R,\xi))=0,\\ \\
\div_\xi W(R,\xi)=0,
\end{array}\right.
\end{equation}
and
\begin{eqnarray}
	\langle W\rangle &=&0,\label{e:MikadoW}\\
    \langle W\otimes W\rangle &=&R.\label{e:MikadoWW}
\end{eqnarray}
\end{lemma}

The first step in the proof of Lemma \ref{l:Mikado} is the following global version of Lemma \ref{l:decomp1} from 
\cite{Nash1954} (that is used to obtain the global decomposition \eqref{e:decomposition}):
\begin{lemma}\label{L:geomNash}
For any compact subset $\mathcal N$ of positive definite $3\times 3$ matrices there exists $\lambda_0\geq1$ and smooth functions $\Gamma_k\in C^\infty(\mathcal N;[0,1])$ for any $k\in\Z^3$ with $|k|\leq\lambda_0$ such that 
\begin{equation}\label{E:Rrepr}
R=\underset{k\in\Z^3,|k|\leq\lambda_0}{\sum}\Gamma_k^2(R)k\otimes k\quad \textrm{ for all }R\in\mathcal{N}.
\end{equation}
\end{lemma}

The proof of Lemma \ref{l:Mikado} is rather simple. 
The vector field $W(R,\cdot)$ will take the form
\begin{equation}\label{E:wdense}
W(R,\xi)=\underset{k\in\Z^3,|k|\leq\lambda_0}{\sum}\Gamma_k(R)\psi_k(\xi)k\, .
\end{equation}
The functions $\psi_k$ are defined as $\psi_k(\xi)=g_k(\mathrm{dist}(\xi, \ell_{k,{p_k}}))$ for some $g_k\in C^\infty_c([0,r_k))$, $r_k>0$, and $\ell_{k,{p_k}}$ is the $\T^3$-periodic extension of the line $\{p_k+tk:\,t\in\R\}$ passing through $p_k$ in direction $k$. Since there are only a finite number of such lines, we may choose $p_k$ and $r_k>0$ in such a way that 
\begin{equation}\label{e:disjoint}
\mathrm{supp }\,\psi_i\cap\mathrm{supp }\,\psi_j=\emptyset\qquad\text{for all $i\neq j$.}
\end{equation}
Thus $W$ consists of a finite collection of disjoint straight tubes such that in each tube $W$ is a straight pipe flow and outside the tubes $W=0$.
In particular $W$ satisfies the stationary ``pressureless'' Euler equations \eqref{e:Mikado}.
Furthermore, the profile functions $g_k$ can be chosen so that $\int_{\T^3}\psi_k(\xi)\,d\xi=0$ and
\[
\fint_{\T^3} \psi_k^2(\xi)\,d\xi=1\quad\textrm{ for all $k$}.
\]
Then \eqref{e:MikadoW} is easily satisfied, and because of \eqref{e:disjoint} we also have
\[
\fint_{\T^3}W\otimes W\,d\xi=\sum_{k}\fint_{\T^3}\Gamma_k^2(R)\psi_k^2(\xi)k\otimes k\,d\xi=\sum_{k}\Gamma_k^2(R)k\otimes k=R.
\]
Therefore \eqref{e:MikadoWW} is satisfied.

This set of flows can be used easily to obtain one initial perturbation of an arbitrary starting subsolution $(v_0,p_0,R_0)$. Indeed, 
the transport ansatz from \eqref{e:Phil} can be used without a temporal cutoff: this time we are not interested in precise estimates for the perturbation, the goal is just to obtain a sufficiently small new Reynolds term $\mathring{R}_1$ so that \eqref{e:Rq} is satisfied. After this single step we then obtain $(v_1,p_1,R_1)$ to which the iteration with Beltrami flows described in the previous sections can be applied. Ironically, at the moment we are not able to carry out the iteration using Mikado flows - the difficulty is in controlling the interaction of two sets of Mikado flows in the temporal overlap regions $I_j\cap I_{j+1}$ (cf.~\eqref{e:overlapslow}).

\section{Further considerations and open questions}\label{s:speculations}

\subsection{Borisov's rigidiy theorem and the threshold $\frac{1}{2}$} In \cite{Pogorelov73} Pogorelov introduced the notion of bounded extrinsic curvature for surfaces in $\mathbb R^3$. Loosely speaking an immersed surface has bounded extrinsic curvature if the area distortion of its Gauss map $N$ is bounded. If the immersion is smooth, this would be a consequence of Gauss' classical theorem, however the definition makes sense as soon as $N$ is a well defined map and thus, for instance, if the immersion is merely $C^1$. 
A consequence of a fundamental result of Pogorelov is the following

\begin{theorem} If $u$ is a $C^1$ immersion of a $2$-dimensional Riemannian manifold $(M,g)$ with positive Gauss curvature and
$u (M)$ has bounded extrinsic curvature in the sense of Pogorelov, then locally the immersed surface is convex.
\end{theorem}

Higher regularity for the immersed surface then follows from the (nowadays classical) regularity theory for the Monge - Amp\`ere equations (cf. \cite{Pogorelov73,Sabitov}). The main point in Borisov's works \cite{Borisov58-1,Borisov58-2,Borisov58-3,Borisov58-4} is to establish that $C^{1, \frac{2}{3}+\varepsilon}$ immersions of surfaces with positive Gauss curvature have bounded extrinsic curvature. 

In \cite{ConDS} Sergio Conti and the two authors observed that Borisov's statement could be recovered from the validity of the following integral identity
\begin{equation}\label{e:ChangeOfVar}
\int_V f(N(x)) \kappa (x)\, dA(x)\;=\;
\int_{\S^2}f(y)\deg(y,V,N)\, d\sigma(y)
\end{equation}
where
\begin{itemize}
\item $V$ is an arbitrary open subset of $M$;
\item $f$ is any bounded function on ${\mathbb S}^2$;
\item $\deg (y, V, N)$ is the Brouwer
degree of the map $N|_V$ at $y$;
\item $dA$ denotes the Riemannian volume form on $(M,g)$;
\item $d\sigma$ is the standard surface measure on ${\mathbb S}^2$.
\end{itemize}
For smooth immersions $u$ \eqref{e:ChangeOfVar} is equivalent to Gauss' theorem. The main point of \cite{ConDS} is that the validity of
\eqref{e:ChangeOfVar} can be extended with little effort to $C^{1,\frac{2}{3}+\varepsilon}$ immersions $u$: if we regularize $u$ by a standard mollification procedure, although a naive computation seems to require $C^{1,1}$ regularity for the convergence of the left hand side, the commutator estimate of Lemma \ref{l:mollify} allows to lower the regularity to $C^{1, \frac{2}{3}+\varepsilon}$. 
We also refer to \cite{BehrOlb} for a partial generalization to hypersurfaces of higher dimension. 

\medskip

There are a number of reasons to believe that this point of view might lower the rigidity threshold to $\frac{1}{2}$. 

First of all if $C\subset \mathbb R^2$ is a $1$-dimensional set and $N: \mathbb R^2 \to \mathbb R^2$ is a $C^{0, \frac{1}{2}+\varepsilon}$ map, then the image $N (C)$ has zero Lebesgue measure. This is a simple elementary fact, cf. \cite{ConDS}.
Moreover, for every bounded open set $U\subset \mathbb R^2$ with Lipschitz boundary, $\deg (\cdot, V, N)\in L^1 (\mathbb R^2)$. This has been proved recently (and independently) by Olbermann \cite{Olbermann} and Z\"ust \cite{Zust}, with rather different arguments. In fact both references have much more general results, valid in several dimensions and general targets: \cite{Olbermann} contains a suitable $L^p$ estimate, whereas, although the arguments  in \cite{Zust} yield only $L^1$ estimates, they allow for different H\"older exponents  for the components of the map. 

In particular, the $C^{1, \frac{1}{2}+\varepsilon}$ regularity  is enough to make sense of the right hand side of \eqref{e:ChangeOfVar} when $V$ has a Lipschitz boundary and $f$ is an arbitrary bounded function: for a general $C^1$ immersion one must instead require that $f$ is compactly supported in $\mathbb S^2\setminus N (\partial V)$. 

Moreover in \cite{Zust} the author has observed that the $L^1$ bound on the degree combined with the computations of \cite{ConDS} is enough to show the following 

\begin{proposition}\label{p:weak}
If $u: M \to \mathbb R^3$ is a $C^{1, \frac{1}{2}+\varepsilon}$ immersion of a smooth $2$-dimensional Riemannian manifold $(M,g)$ and $f=1$, then the identity \eqref{e:ChangeOfVar} is valid for any open subset $V\subset M$ with Lipschitz boundary. 
\end{proposition}

The rigidity threshold could then be lowered to $\frac{1}{2}$ if the following conjecture were true (Z\"ust in \cite{Zust} has proposed an argument for the conjecture, but unfortunately it contains a crucial gap).

\begin{conjecture}
Assume $N: \mathbb R^2 \supset \Omega \to \mathbb R^2$ is map in $C^{\frac{1}{2}+\eps}$ with the property that
\[
\int \deg (y, N, V)\, dy \geq 0
\] 
for every open $V \subset \subset \Omega$ with Lipschitz boundary. Then $\deg (y, N, V)$ is nonnegative for every open $V\subset\subset  \Omega$ and every $y\not\in N (\partial V)$. 
\end{conjecture}

\subsection{Further results on incompressible Euler and other equations} The techniques introduced in the papers \cite{DS-Inv,DS-JEMS} have been extended to prove several other results in the incompressible Euler equations and for other equations in fluid dynamics. 

Concerning the Euler equations, in \cite{Choffrut} Choffrut showed that the same tools can be suitably modified to produce H\"older continuous dissipative weak solutions when the space domain is the $2$-dimensional torus $\mathbb T^2$.  In \cite{Daneri} Daneri gave a first construction which produces infinitely many solutions with the same initial data and have nonincreasing energy. This result was improved further in \cite{DaneriSz}, where the H\"older regularity of \cite{Daneri} has been pushed to match that of Theorem \ref{t:onsager}(ii). The same paper also shows that the initial data allowing for such nonuniqueness theorem are indeed dense in $L^2$.  In \cite{IsettOh} Isett and Oh produced H\"older solutions which are compactly supported in time and {\em space} when the space domain is $\mathbb R^3$. 

Remarkably, in \cite{IsettVicol} Isett and Vicol have succeeded in implementing a multistep iteration scheme which produces H\"older continuous solutions to active scalar equations when the multiplier is not odd. This combines the ideas of \cite{DS-Inv} with previous techniques used in \cite{CFG,Sz-IPM,Shvydkoy} to produce bounded solutions when the multiplier is even. In \cite{Tao1,Tao2} Tao and Zhang have extended the 
techniques of \cite{DS-Inv,DS-JEMS} to produce similar results for the Boussinesq Equation. 

In the cases of bounded solutions it has been shown in \cite{DS-ARMA} that convex integration can be used to produce very irregular solutions which satisfy the local energy inequality
\begin{equation}\label{e:local}
\partial_t \frac{|u|^2}{2} + {\rm div}\, \left(\left(\frac{|u|^2}{2} + p\right) u\right) \leq 0\, 
\end{equation}
and that therefore the latter condition is still not enough to ensure uniqueness of a weak solution. This remains true even for initial data which have very mild discontinuities, as shown by the second author in \cite{Sz-IPM}. In fact the same constructions can be used in compressible fluid dynamics to disprove the uniqueness of entropy admissible weak solutions for some regular (more precisely Lipschitz) initial data, cf. \cite{DS-ARMA,CDK,ChK}. It is presently not known whether one could use techniques similar to those of \cite{DS-Inv} to construct {\em continuous} solutions which satisfy the local energy inequality \eqref{e:local}. In particular it would be of some interest
to disprove the uniqueness of {\em piecewise continuous} entropy admissible weak solutions in compressible fluid dynamics.

\bibliographystyle{acm}
\bibliography{Nash-Bulletin}

\end{document}